\documentclass{article}%
\usepackage{amssymb}
\usepackage{amsmath}
\usepackage{amsfonts}
\usepackage{graphicx}%
\providecommand{\U}[1]{\protect\rule{.1in}{.1in}}
\newtheorem{theorem}{Theorem}

\newtheorem{definition}[theorem]{Definition}

\newtheorem{lemma}[theorem]{Lemma}

\newtheorem{proposition}[theorem]{Proposition}
\newtheorem{remark}[theorem]{Remark}

\begin{document}

\title{A minimization method and applications to the study of solitons}
\author{Vieri Benci\thanks{ Dipartimento di Matematica Applicata, Universit\`a degli
Studi di Pisa, Via F. Buonarroti 1/c, Pisa, ITALY and Department of
Mathematics, College of Science, King Saud University, Riyadh, 11451, SAUDI
ARABIA. e-mail: \texttt{benci@dma.unipi.it}}
\and Donato Fortunato\thanks{Dipartimento di Matematica, Universit\`{a} degli Studi
di Bari Aldo Moro, Via Orabona 4, 70125 Bari, Italy,
e-mail:\texttt{\ fortunat@dm.uniba.it}}
\and \bigskip
\and \textit{Dedicated to V. Lakshmikantham}}
\maketitle

\begin{abstract}
Roughly speaking a \textit{solitary wave} is a solution of a field equation
whose energy travels as a localized packet and which preserves this
localization in time. A \textit{soliton} is a solitary wave which exhibits
some strong form of stability so that it has a particle-like behavior. In this
paper, we prove a general, abstract theorem (Theorem \ref{astra1}) which
allows to prove the existence of a class of solitons. Such solitons are
suitable minimizers of a constrained functional and they are called
hylomorphic solitons. Then we apply the abstract theory to problems related to
the nonlinear Schr\"{o}dinger equation (NSE) and to the nonlinear Klein-Gordon
equation (NKG).

\bigskip


\bigskip


\end{abstract}
\tableofcontents

\bigskip

\section{Introduction}

In some recent papers (\cite{BBBM},\ \cite{BBGM}, \cite{milano},
\cite{befobeam}, \cite{befolat}, \cite{befoQ}) the existence of solitons has
been proved using variational methods. In this paper, we prove a general,
abstract theorem (Theorem \ref{astra1}) which applies to most of the
situations analyzed in the mentioned papers.

The proof of Theorem \ref{astra1} is carried out in two steps: in the first
step (section \ref{mr}) the research of the minimizers of a constrained
functional is reduced to the study of the minimizers of a suitable free
functional. The existence of such minimizers is stated in Theorem \ref{furbo}.
In the second step stability properties of these minimizers are proved
(section \ref{ehs}).

These two theorems (Theorem \ref{furbo} and Theorem \ref{astra1}) can be
applied to the situations described in the quoted papers relative to solitons.
In section \ref{ehs}, we give an abstract definition of \textit{soliton}
(Definition \ref{ds}) and \textit{hylomorphic soliton }(Definition \ref{tdc}).
These solitons are stable minimizers of a constrained functional. Then we
apply the abstract theorems \ref{astra1} and \ref{furbo} to new problems
related to the nonlinear Schr\"{o}dinger equation (NSE) and to the nonlinear
Klein-Gordon equation (NKG): in section \ref{EV} we use Theorem \ref{furbo} in
order to prove the existence of vortices for the NSE with a potential $V$
which is periodic in one direction (Theorem \ref{vortices}). In section
\ref{periodic} we use Theorem \ref{astra1} in order to prove the existence of
hylomorphic solitons for NSE in a lattice, namely in presence of a periodic
potential $V$.

Finally in section \ref{kleingordon} we use Theorem \ref{astra1} to prove the
existence of hylomorphic solitons for the nonlinear Klein-Gordon equation (NKG).

\section{A minimization result\label{mr}}

Let $E$ and $C$ be two functionals on an Hilbert space $X.$ We are interested
in the following minimization problem: find values of $\sigma\in\mathbb{R}$,
such that $E$ attains a constrained minimum on $\mathfrak{M}_{\sigma}$ where
\[
\mathfrak{M}_{\sigma}:=\left\{  \mathbf{u}\in X\ |\ \left\vert C(\mathbf{u}%
)\right\vert =\sigma\right\}  .
\]
In the next section, we will describe the abstract framework where we will
work and then we will prove an abstract existence theorem.

\subsection{The abstract framework\label{af}}

We need some definitions. These definitions are related to a couple $(X,G)$
where $G$ is a group acting on the Hilbert space $X.$ In our applications $G$
will be a subgroup of the group of translations.

\begin{definition}
A non-empty subset $\Gamma\subset X$ is called $G$-invariant if
\[
\forall\mathbf{u}\in\Gamma,\ \forall g\in G,\ g\mathbf{u}\in\Gamma.
\]

\end{definition}

\begin{definition}
A functional $J$ on $X$ is called $G$\emph{-invariant} if
\[
\forall g\in G,\text{ }\forall\mathbf{u}\in X,\ J\left(  g\mathbf{u}\right)
=J\left(  \mathbf{u}\right)  .
\]

\end{definition}

\begin{definition}
\label{gcompatto}A closed $G$-invariant set $\Gamma\subset X$ is called
$G$-compact if for any sequence $\mathbf{u}_{n}$ in $\Gamma$ there is a
sequence $g_{n}\in G,$ such that $g_{n}\mathbf{u}_{n}$ has a converging
subsequence. Clearly a sequence $\mathbf{u}_{n}$ in $X$ will be called
$G$\emph{-compact }if its image is $G$-compact.
\end{definition}

\begin{definition}
\label{gcompattoa}A $G$-invariant functional $J$ on $X$ is called $G-$compact
if any minimizing sequence $\mathbf{u}_{n}$ is $G-$compact.
\end{definition}

Clearly a $G$-compact functional has a $G$-compact set of minimizers.

In order to prove an existence result for the minimizers of $E$ on
$\mathfrak{M}_{\sigma}$, we impose some assumptions to $E$ and $C$; to do this
we need some other definitions:

\begin{definition}
\label{splittinga}(\textbf{Splitting property)}We say that a functional $F$ on
$X$ has the splitting property if given a sequence $\mathbf{u}_{n}%
=\mathbf{u}+\mathbf{w}_{n}$ in $X$ such that $\mathbf{w}_{n}$ converges weakly
to $0$, we have that
\begin{equation}
F(\mathbf{u}_{n})=F(\mathbf{u})+F(\mathbf{w}_{n})+o(1).
\end{equation}

\end{definition}

\begin{remark}
\label{quadratic}\textbf{\ }A symmetric, continuous quadratic form satisfies
the splitting property; in fact, in this case, we have that $F(\mathbf{u}%
):=\left\langle L\mathbf{u},\mathbf{u}\right\rangle $ for some continuous
selfajoint operator $L;$ then given a sequence $\mathbf{u}_{n}=\mathbf{u}%
+\mathbf{w}_{n}$ with $\mathbf{w}_{n}\rightharpoonup0$ weakly, we have that%
\begin{align*}
F(\mathbf{u}_{n})  &  =\left\langle L\mathbf{u},\mathbf{u}\right\rangle
+\left\langle L\mathbf{w}_{n},\mathbf{w}_{n}\right\rangle +2\left\langle
L\mathbf{u},\mathbf{w}_{n}\right\rangle \\
&  =F(\mathbf{u})+F(\mathbf{w}_{n})+o(1).
\end{align*}

\end{remark}

Now we can formulate the required properties on $E$ and $C$:

\begin{itemize}
\item \textit{(EC-1) \textbf{(Value at 0)} }$E,$ $C$ \textit{are} $C^{1}$,
bounded functionals such that\textit{\ }%
\[
E(0)=0,\ C(0)=0;\ E^{\prime}(0)=0;\ C^{\prime}(0)=0.
\]

\item \textit{(EC-2) \textbf{(Invariance)\ }}$E$ \textit{and }$C$\textit{\ are
}$G$\textit{-invariant.}

\item \textit{(EC-3)\textbf{(Coercivity)} We distinguish two cases: } $C\geq0$
\textit{and }$C\ $\textit{not positive. If }$C\geq0$\textit{\ we assume that}
\textit{there exists }$a\geq0$ \textit{and }$s\geq1$\textit{\ such that}

\begin{itemize}
\item (i) $E(\mathbf{u})+aC(\mathbf{u})^{s}\geq0;$

\item (ii) \textit{if }$\left\Vert \mathbf{u}\right\Vert \rightarrow\infty
,\ $\textit{then} $E(\mathbf{u})+aC(\mathbf{u})^{s}\rightarrow\infty;$

\item (iii) \textit{for any} \textit{bounded sequence }$\mathbf{u}_{n}$
\textit{in }$X$ \textit{such that} $E(\mathbf{u}_{n})+aC(\mathbf{u}_{n}%
)^{s}\rightarrow0,\ $\textit{we have that }$\mathbf{u}_{n}\rightarrow0. $
\end{itemize}

\textit{In the case in which }$C$\textit{\ is not positive we assume that (i),
(ii), (iii) hold true with }$a=0.$

\item \textit{(EC-4)\textbf{(Splitting property)} }$E$\textit{\ and }%
$C$\textit{\ satisfy the splitting property.}
\end{itemize}

Before stating our main theorem we need this definition:

\begin{definition}
\label{na}\textbf{\ }A bounded sequence $\mathbf{u}_{n}$ in $X$ is called
\textbf{vanishing sequence}\textit{\ if for any subsequence }$\mathbf{u}%
_{n_{k}}$ and any sequence $g_{k}\subset G$, $g_{k}\mathbf{u}_{n_{k}}$
converges weakly to $0.$
\end{definition}

Observe that the notion of vanishing sequence depends on the group $G$ acting
on $X.$ Clearly a bounded sequence $\mathbf{u}_{n}$ in $X$ is a non-vanishing
sequence\textit{\ }if there exists a subsequence\textit{\ }$\mathbf{u}_{n_{k}%
}$ and a sequence $g_{k}\subset G$, such that  $g_{k}\mathbf{u}_{n_{k}}$
converges weakly to some $\mathbf{u}\neq0.$

So, if $\mathbf{u}_{n}\rightarrow0$ strongly, then $\mathbf{u}_{n}$ is a
vanishing sequence. However, if $\mathbf{u}_{n}\rightharpoonup0$ weakly, it
might happen that it is a non-vanishing sequence.

\subsection{The minimization theorem}

We now set
\[
\Lambda(\mathbf{u}):=\frac{E(\mathbf{u})}{\left\vert C(\mathbf{u})\right\vert
}%
\]
and
\begin{equation}
\Lambda_{0}:=\ \inf\left\{  \lim\inf\ \Lambda(\mathbf{u}_{n})\ |\ \mathbf{u}%
_{n}\ \text{is a vanishing sequence}\right\}  . \label{hylo}%
\end{equation}
Now we state and prove the following existence result:

\begin{theorem}
\label{furbo} Let $E$ and $C$ be two functionals on a Hilbert space $X$ and
$G$ be a group acting on $X.$ Assume that $E$ and $C$ satisfy
(EC-1),...,(EC-4) and
\begin{equation}
\underset{}{\inf}\Lambda(\mathbf{u})<\Lambda_{0}. \label{hh}%
\end{equation}
Then there is $\bar{\delta}>0$ and a family of values $c_{\delta},\ \delta
\in\left(  0,\bar{\delta}\right)  ,$ such that the minimum
\begin{equation}
e_{\delta}=\min\left\{  E(\mathbf{u})\ |\ \left\vert C(\mathbf{u})\right\vert
=c_{\delta}\right\}
\end{equation}
exists and the set $\Gamma_{c_{\delta}}$ of minimizers is $G$-compact.
Moreover $\Gamma_{c_{\delta}}\ $can be characterized as the set of minimizers
of the functional%
\[
J_{\delta}(\mathbf{u})=\Lambda\left(  \mathbf{u}\right)  +\delta
\Phi(\mathbf{u})\text{ }\delta\in\left(  0,\bar{\delta}\right)
\]
where $\Phi(\mathbf{u})=E(\mathbf{u})+2a\left\vert C(\mathbf{u})\right\vert
^{s}$ with $a$ as in (EC-3).
\end{theorem}

\begin{remark}
\label{valutazione}When we will apply Th. \ref{furbo}, it is necessary to
estimate $\Lambda_{0};$ the following inequalities may help to do this. In
order to give an estimate from below, assume that there exists a seminorm
$\left\Vert \mathbf{u}\right\Vert _{\sharp}$ such that
\begin{equation}
\left\{  \mathbf{u}_{n}\ \text{is a vanishing sequence}\right\}
\Rightarrow\left\Vert \mathbf{u}_{n}\right\Vert _{\sharp}\rightarrow0.
\label{seminorm}%
\end{equation}
Then we have that%
\begin{equation}
\Lambda_{0}\geq\ \underset{\left\Vert \mathbf{u}\right\Vert _{\sharp
}\rightarrow0}{\lim\inf}\ \Lambda(\mathbf{u}). \label{aaaa}%
\end{equation}
Let $\left\Vert {}\right\Vert $ denote the norm in $X.$ Since
\[
\left\Vert \mathbf{u}_{n}\right\Vert \rightarrow0\Rightarrow\left\{
\mathbf{u}_{n}\ \text{is a vanishing sequence}\right\}  ,
\]
we have that%
\begin{equation}
\Lambda_{0}\leq\ \underset{\left\Vert \mathbf{u}\right\Vert \rightarrow0}%
{\lim\inf}\ \Lambda(\mathbf{u}). \label{bbbb}%
\end{equation}
Moreover, if $E$ and $C$ are twice differentiable in $0$, by (EC-0), we have
that%
\begin{align*}
\underset{\left\Vert \mathbf{u}\right\Vert \rightarrow0}{\lim\inf}%
\ \Lambda(\mathbf{u})  &  =\ \underset{\left\Vert \mathbf{u}\right\Vert
\rightarrow0}{\lim\inf}\frac{E(0)+E^{\prime}(0)\left[  \mathbf{u}\right]
+E^{\prime\prime}(0)\left[  \mathbf{u,u}\right]  +o(\left\Vert \mathbf{u}%
\right\Vert ^{2})}{\left\vert C(0)+C^{\prime}(0)\left[  \mathbf{u}\right]
+C^{\prime\prime}(0)\left[  \mathbf{u,u}\right]  +o(\left\Vert \mathbf{u}%
\right\Vert ^{2})\right\vert }\\
&  =\ \underset{}{\inf}\frac{E^{\prime\prime}(0)\left[  \mathbf{u,u}\right]
}{\left\vert C^{\prime\prime}(0)\left[  \mathbf{u,u}\right]  \right\vert }%
\end{align*}
and hence%
\[
\Lambda_{0}\leq\ \underset{}{\inf}\frac{E^{\prime\prime}(0)\left[
\mathbf{u,u}\right]  }{\left\vert C^{\prime\prime}(0)\left[  \mathbf{u,u}%
\right]  \right\vert }.
\]
It is can be seen that in many applications the two limits (\ref{aaaa}) and
(\ref{bbbb}) coincide and in this case we get a sharp estimate for
$\Lambda_{0}.$
\end{remark}

\bigskip

\begin{remark}
If $X=H^{1}(\mathbb{R}^{N})$ and $G=\mathbb{Z}^{N}$ with the action
(\ref{group}), then (see Proposition \ref{nonvanishing2}) the norm $\left\Vert
\mathbf{u}\right\Vert _{\sharp}=\left\Vert \mathbf{u}\right\Vert _{L^{t}}$,
$t\in\left(  2,\frac{2N}{N-2}\right)  ,N\geq3,$ satisfies the property
(\ref{seminorm}) of the preceding remark \ref{valutazione} and consequently it
satisfies also (\ref{aaaa})).
\end{remark}

\begin{remark}
The fact that the minimization problem of $E$ on $\mathfrak{M}_{\sigma}$
reduces to the minimization of a free functional $J_{\delta}(\mathbf{u})\ $is
very useful in numerical simulation.
\end{remark}

\bigskip

\textbf{Proof }of theorem \ref{furbo}. By (\ref{hh}) there exists
$\mathbf{v}\in X$ such that $\Lambda(\mathbf{v})<\Lambda_{0},$ then we can
take $\delta>0$ so small that
\[
J_{\delta}(\mathbf{v})=\Lambda\left(  \mathbf{v}\right)  +\delta
\Phi(\mathbf{v})<\Lambda_{0}\text{ }%
\]
and define
\[
\bar{\delta}=\sup\left\{  \delta\ |\ \exists\mathbf{v}:\Lambda\left(
\mathbf{v}\right)  +\delta\Phi(\mathbf{v})<\Lambda_{0}\ \right\}  .
\]

Then
\begin{equation}
\underset{\mathbf{u}\in X}{\inf}J_{\delta}\left(  \mathbf{u}\right)
<\Lambda_{0}\text{ for }\delta\in\left(  0,\bar{\delta}\right)  .
\label{favino}%
\end{equation}

Now we show that
\begin{equation}
J_{\delta}\left(  \mathbf{u}\right)  \geq\frac{\delta}{2}\Phi(\mathbf{u})-M
\label{pli}%
\end{equation}
where $M$ is a suitable constant. Clearly, if $C$ is not positive, we have
$a=0$ in (EC-3)(i); then $\Phi(\mathbf{u})=E(\mathbf{u})\geq0$ and $J_{\delta
}(\mathbf{u})=\frac{E(\mathbf{u})}{\left\vert C(\mathbf{u})\right\vert
}+\delta E(\mathbf{u}).$ Then (\ref{pli}) is obviously satisfied.

Now assume that $C(\mathbf{u})\geq0.$ By (EC-3)(i) we have that%
\begin{equation}
E(\mathbf{u})\geq-aC(\mathbf{u})^{s} \label{aaa}%
\end{equation}
and hence%
\begin{equation}
\frac{E(\mathbf{u})}{C(\mathbf{u})}\geq-aC(\mathbf{u})^{s-1}. \label{bibi}%
\end{equation}

Then, by (\ref{aaa}) and (\ref{bibi}), we get
\begin{align*}
J_{\delta}(\mathbf{u})  &  =\frac{E(\mathbf{u})}{C(\mathbf{u})}+\delta
\Phi(\mathbf{u})\geq-aC(\mathbf{u})^{s-1}+\frac{\delta}{2}\left[
E(\mathbf{u})+2aC(\mathbf{u})^{s}\right]  +\frac{\delta}{2}\Phi(\mathbf{u})\\
&  \geq-aC(\mathbf{u})^{s-1}+\frac{\delta}{2}\left[  -aC(\mathbf{u}%
)^{s}+2aC(\mathbf{u})^{s}\right]  +\frac{\delta}{2}\Phi(\mathbf{u})\\
&  \geq-aC(\mathbf{u})^{s-1}+\frac{a\delta}{2}C(\mathbf{u})^{s}+\frac{\delta
}{2}\Phi(\mathbf{u})\geq\frac{\delta}{2}\Phi(\mathbf{u})-M
\end{align*}
where%
\[
M=\ -a\underset{t\geq0}{\min}\left(  \frac{\delta}{2}t^{s}-t^{s-1}\right)  .
\]
Then (\ref{pli}) has been proved.

Now let us prove that $J_{\delta}$ is $G-$ compact (see Definition
\ref{gcompattoa}).

Let $\mathbf{u}_{n}$ be a minimizing sequence of $J_{\delta}.$ This sequence
$\mathbf{u}_{n}$ is bounded in $X.$ In fact, arguing by contradiction, assume
that, up to a subsequence, $\left\Vert \mathbf{u}_{n}\right\Vert $
$\rightarrow\infty.$ Then by (\ref{pli}) and (EC-3) (ii) we get $J_{\delta
}(\mathbf{u}_{n})\rightarrow\infty$ which contradicts the fact that
$\mathbf{u}_{n}$ is a minimizing sequence of $J_{\delta}.$

Since $\mathbf{u}_{n}$ is minimizing for $J_{\delta},$ by (\ref{favino}),
there exists $\eta>0$ such that, for $n$ sufficiently large,
\[
\frac{E(\mathbf{u}_{n})}{\left\vert C(\mathbf{u}_{n})\right\vert }+\delta
\Phi(\mathbf{u}_{n})<\Lambda_{0}-\eta
\]
then
\begin{equation}
\Lambda(\mathbf{u}_{n})=\frac{E(\mathbf{u}_{n})}{\left\vert C(\mathbf{u}%
_{n})\right\vert }<\Lambda_{0}-\eta. \label{upper}%
\end{equation}

On the other hand
\begin{equation}
\Lambda(\mathbf{u}_{n})\text{ is bounded below.} \label{below}%
\end{equation}

In fact: since $\Phi$ is bounded and $\mathbf{u}_{n}$is bounded in $X$,
$\Phi(\mathbf{u}_{n})$ is bounded. So we deduce from (\ref{pli}) that
$J_{\delta}(\mathbf{u}_{n})$ is bounded below. Then $\Lambda(\mathbf{u}%
_{n})=J_{\delta}(\mathbf{u}_{n})-\delta\Phi(\mathbf{u}_{n})$ is bounded below
and so (\ref{below}) is proved. By (\ref{upper}) and (\ref{below}) we have,
for some subsequence, that
\begin{equation}
\Lambda(\mathbf{u}_{n})\rightarrow\lambda,\text{ }-\infty<\lambda<\Lambda_{0}.
\label{bbb}%
\end{equation}
Then, by (\ref{hylo}), $\mathbf{u}_{n}$ is a bounded non vanishing sequence.
Hence, by Def. \ref{na}, we can extract a subsequence $\mathbf{u}_{n_{k}}$ and
we can take a sequence $g_{k}\subset G$ such that $\mathbf{u}_{k}^{\prime
}:=g_{k}\mathbf{u}_{n_{k}}$ is weakly convergent to some
\begin{equation}
\mathbf{\bar{u}}\neq0. \label{agg}%
\end{equation}
We can write
\[
\mathbf{u}_{n}^{\prime}=\mathbf{\bar{u}}+\mathbf{w}_{n}%
\]
with $\mathbf{w}_{n}\rightharpoonup0$ weakly. We want to prove that
$\mathbf{w}_{n}\rightarrow0$ strongly. First of all we will show that%
\begin{equation}
\lim\Phi(\mathbf{\bar{u}}+\mathbf{w}_{n})\geq\Phi(\mathbf{\bar{u}})+\lim
\Phi(\mathbf{w}_{n}). \label{good}%
\end{equation}

If $C(\mathbf{u})$ is not positive we have $a=0$ in (EC-3)(i), then
$\Phi(\mathbf{u})=E(\mathbf{u})$ and clearly (\ref{good}) holds as an equality
since $E$ satisfies assumption (EC-4) (splitting property).

Now assume that $C(\mathbf{u})\geq0.$ Then by (EC-4) and since $s\geq1,$ we
have that
\begin{align}
\lim\Phi(\mathbf{\bar{u}}+\mathbf{w}_{n})  &  =\lim\left(  E(\mathbf{\bar{u}%
}+\mathbf{w}_{n})+2aC(\mathbf{\bar{u}}+\mathbf{w}_{n})^{s}\right) \nonumber\\
&  =E(\mathbf{\bar{u}})+\lim E(\mathbf{w}_{n})+2a\lim\left(  C(\mathbf{\bar
{u}})+C(\mathbf{w}_{n})\right)  ^{s}\nonumber\\
&  \geq E(\mathbf{\bar{u}})+\lim E(\mathbf{w}_{n})+2a\lim\left(
C(\mathbf{\bar{u}})^{s}+C(\mathbf{w}_{n})^{s}\right) \nonumber\\
&  =E(\mathbf{\bar{u}})+2aC(\mathbf{\bar{u}})^{s}+\lim E(\mathbf{w}%
_{n})+2a\lim C(\mathbf{w}_{n})^{s}\nonumber\\
&  =\Phi(\mathbf{\bar{u}})+\lim\Phi(\mathbf{w}_{n}). \label{qw}%
\end{align}
So (\ref{good}) has been proved.

Next we show that%
\begin{equation}
C(\mathbf{\bar{u}+\mathbf{w}_{n})}\text{ does not converge to }0. \label{not}%
\end{equation}
Arguing by contradiction assume that $C(\mathbf{\bar{u}+\mathbf{w}_{n})}$
converges to $0.$ Then, since $\mathbf{\bar{u}}+\mathbf{\mathbf{w}_{n}}$ is a
minimizing sequence for $J_{\delta},$ also $E(\mathbf{\bar{u}+\mathbf{w}_{n}%
)}$ converges to $0$ and then
\[
E(\mathbf{\bar{u}+\mathbf{w}_{n})+}a\left\vert C(\mathbf{\bar{u}%
+\mathbf{w}_{n})}\right\vert ^{s}\rightarrow0.
\]

So, by (EC-3)(iii), we get
\begin{equation}
\mathbf{\bar{u}+\mathbf{w}_{n}}\rightarrow0\text{ in }X. \label{s}%
\end{equation}
From (\ref{s}) and since $\mathbf{\mathbf{w}_{n}}\rightharpoonup0$ weakly in
$X,$ we have that $\mathbf{\bar{u}}=0$, contradicting (\ref{agg}). Then
(\ref{not}) holds and consequently, up to a subsequence, we have
\begin{equation}
\left\vert C(\mathbf{\bar{u}+\mathbf{w}_{n})}\right\vert =\left\vert
C(\mathbf{\bar{u}})+C(\mathbf{w}_{n})+o(1)\right\vert \geq const.>0.
\label{fir}%
\end{equation}

Now, we set%
\begin{align*}
j_{\delta}  &  =\inf J_{\delta}=\lim J_{\delta}(\mathbf{u}_{n}^{\prime
});\ e_{\delta}=E(\mathbf{\bar{u}});\ c_{\delta}=\left\vert C(\mathbf{\bar{u}%
})\right\vert \\
e_{1}  &  =\lim E(\mathbf{w}_{n});\ c_{1}=\lim\left\vert C(\mathbf{w}%
_{n})\right\vert .
\end{align*}
Observe that the limits $\lim E(\mathbf{w}_{n})$ and $\lim\left\vert
C(\mathbf{w}_{n})\right\vert $ exist (up to subsequences), since $E$ and $C$
are bounded functionals and $\mathbf{w}_{n}$ weakly converges.

Now we have
\begin{equation}
\lim\frac{E(\mathbf{\bar{u}})+E(\mathbf{w}_{n})+o(1)}{\left\vert
C(\mathbf{\bar{u}})+C(\mathbf{w}_{n})+o(1)\right\vert }\geq\frac{\ e_{\delta
}+e_{1}}{c_{\delta}+c_{1}}. \label{anche}%
\end{equation}

In fact, as usual we distinguish two cases: if $C\geq0$ (\ref{anche}) holds
since%
\[
\lim\left(  C(\mathbf{\bar{u}})+C(\mathbf{w}_{n})+o(1)\right)  =c_{\delta
}+c_{1}.
\]
On the other hand, if $C$ is not positive, by (EC-3)(i), we have $E\geq0$ and
in this case
\[
\lim\left\vert C(\mathbf{\bar{u}})+C(\mathbf{w}_{n})+o(1)\right\vert \leq
c_{\delta}+c_{1}.
\]
So (\ref{anche}) holds also when $C$ is not positive.

Now by (\ref{anche}) and (\ref{good}), we have that%
\begin{align}
j_{\delta}  &  =\lim\left[  \frac{E(\mathbf{u}_{n}^{\prime})}{\left\vert
C(\mathbf{u}_{n}^{\prime})\right\vert }+\delta\Phi(\mathbf{u}_{n}^{\prime
})\right] \nonumber\\
&  =\lim\frac{E(\mathbf{\bar{u}})+E(\mathbf{w}_{n})+o(1)}{\left\vert
C(\mathbf{\bar{u}})+C(\mathbf{w}_{n})+o(1)\right\vert }+\delta\lim
\Phi(\mathbf{\bar{u}}+\mathbf{w}_{n})\nonumber\\
&  \geq\frac{\ e_{\delta}+e_{1}}{c_{\delta}+c_{1}}+\delta\lim\Phi
(\mathbf{w}_{n})+\delta\Phi(\mathbf{\bar{u}}).\ \label{pipa}%
\end{align}
Now we want to prove that
\begin{equation}
\frac{\ e_{1}}{c_{1}}\geq\frac{\ e_{\delta}}{c_{\delta}}. \label{impo2}%
\end{equation}

We argue indirectly and we suppose that
\begin{equation}
\frac{\ e_{\delta}}{c_{\delta}}>\frac{\ e_{1}}{c_{1}}. \label{impo}%
\end{equation}
By the above inequality it follows that
\begin{equation}
\frac{\ e_{\delta}+e_{1}}{c_{\delta}+c_{1}}=\frac{\frac{e_{\delta}}{c_{\delta
}}\ c_{\delta}+\frac{e_{1}}{c_{1}}\ c_{1}}{c_{\delta}+c_{1}}>\frac{\frac
{e_{1}}{c_{1}}\ c_{\delta}+\frac{e_{1}}{c_{1}}\ c_{1}}{c_{\delta}+c_{1}}%
=\frac{e_{1}}{c_{1}} \label{ava}%
\end{equation}
and hence, by (\ref{pipa}) and (\ref{ava}), we get%
\[
j_{\delta}\geq\frac{\ e_{\delta}+e_{1}}{c_{\delta}+c_{1}}+\delta\lim
\Phi(\mathbf{w}_{n})+\delta\Phi(\mathbf{\bar{u}})\
\]%
\[
>\frac{e_{1}}{c_{1}}+\delta\lim\Phi(\mathbf{w}_{n})+\delta\Phi(\mathbf{\bar
{u}})=\lim J_{\delta}(\mathbf{w}_{n})+\delta\Phi(\mathbf{\bar{u}})
\]%
\[
\geq\inf J_{\delta}+\delta\Phi(\mathbf{\bar{u}})=j_{\delta}+\delta
\Phi(\mathbf{\bar{u}})\geq j_{\delta}.
\]
So we get a contradiction and (\ref{impo}) cannot occur. Then we have
(\ref{impo2}).

Now, by (\ref{impo2}), we get%
\[
\frac{\ e_{\delta}+e_{1}}{c_{\delta}+c_{1}}=\frac{\frac{e_{\delta}}{c_{\delta
}}\ c_{\delta}+\frac{e_{1}}{c_{1}}\ c_{1}}{c_{\delta}+c_{1}}\geq\frac
{\frac{e_{\delta}}{c_{\delta}}\ c_{\delta}+\frac{e_{\delta}}{c_{\delta}%
}\ c_{1}}{c_{\delta}+c_{1}}=\frac{e_{\delta}}{c_{\delta}}.
\]
So we get%

\[
\frac{\ e_{\delta}+e_{1}}{c_{\delta}+c_{1}}\geq\frac{e_{\delta}}{c_{\delta}}.
\]
Then, using (\ref{pipa}), the above inequality and the fact that $j_{\delta
}=\inf J_{\delta},$we get
\begin{align*}
j_{\delta}  &  \geq\frac{\ e_{\delta}+e_{1}}{c_{\delta}+c_{1}}+\delta
\Phi(\mathbf{\bar{u}})+\delta\lim\Phi(\mathbf{w}_{n})\\
&  \geq\frac{e_{\delta}}{c_{\delta}}+\delta\Phi(\mathbf{\bar{u}})+\delta
\lim\Phi(\mathbf{w}_{n})\\
&  =J_{\delta}(\mathbf{\bar{u}})+\delta\lim\Phi(\mathbf{w}_{n})\geq j_{\delta
}+\delta\lim\Phi(\mathbf{w}_{n}).
\end{align*}
Then%
\[
\delta\lim\Phi(\mathbf{w}_{n})\leq0
\]
and, by (EC-3)(iii), $\mathbf{w}_{n}\rightarrow0\ $and hence $\mathbf{u}%
_{n}^{\prime}\rightarrow\mathbf{\bar{u}}$ strongly and $\mathbf{\bar{u}}$ is a
minimizer of $J_{\delta}.$

So we conclude that $J_{\delta}$ is $G-$compact.

Since $\mathbf{\bar{u}}$ is a minimizer of $J_{\delta}$,
clearly$\ \mathbf{\bar{u}}$ minimizes also the functional%
\[
\frac{E(\mathbf{u})}{c_{\delta}}+\delta\left[  E(\mathbf{u})+ac_{\delta}%
^{s}\right]  =\left(  \frac{1}{c_{\delta}}+\delta\right)  E(\mathbf{u})+\delta
ac_{\delta}^{s}%
\]
on the set $\left\{  \mathbf{u}\in X\ |\ \left\vert C(\mathbf{u})\right\vert
=c_{\delta}\right\}  $ and hence $\mathbf{\bar{u}}$ minimizes also $E$ on this
set. Now denote by $\Gamma_{c_{\delta}}$ the set of such minimizers. It is
easy to see that viceversa $\Gamma_{c_{\delta}}$ is contained in the set of
minimizers of $J_{\delta}.$ So, since $J_{\delta}$ is $G$-compact, we conclude
that $\Gamma_{c_{\delta}}$ is $G$-compact.

$\square$

\section{An existence result of vortices for NSE\label{EV}}

The existence of vortices is an interesting and old issue in many questions of
mathematical physics as superconductivity, classical and quantum field theory,
string and elementary particle theory (see the pioneering papers \cite{ab},
\cite{nil} and e.g. the more recent ones \cite{Kim93}, \cite{Volk},
\cite{Volk-Wohn}, \cite{vil} with their references).

From mathematical viewpoint, the existence of vortices for the nonlinear
Klein-Gordon equations (NKG), for nonlinear Schroedinger equations (NSE) and
for nonlinear Klein-Gordon-Maxwell\ equations (NKGM) has been studied in some
recent papers ( \cite{bevi}, \cite{BBR07}, \cite{bbr08}, \cite{befov07},
\cite{befocom}, \cite{bebo}, \cite{bebosi}, \cite{BF09TA}).

Many of the previous results can be obtained applying Th.\ref{furbo}. Here we
will consider a case not covered by the existing literature, namely the study
of vortices in NSE when the potential $V(x)$ depends only on the third
variable and it is periodic, namely for all $k\in\mathbb{Z}$
\begin{equation}
V(x_{1},x_{2},x_{3})=V(x_{3})=V(x_{3}+k). \label{per}%
\end{equation}

\subsection{Statement of the problem\label{beginning}}

Let us consider the nonlinear Schroedinger equation:%
\begin{equation}
i\frac{\partial\psi}{\partial t}=-\frac{1}{2}\Delta\psi+\frac{1}{2}W^{\prime
}\left(  \psi\right)  +V(x)\psi. \tag{NSE}\label{KG}%
\end{equation}
where $\psi(t,x)$ is a complex valued function defined on the space-time
$\mathbb{R\times R}^{N}$ $(N\geq3),\ V:\mathbb{R}^{N}\mathbb{\rightarrow
R},\ W:\mathbb{C\rightarrow R}$ such that $W(\psi)=F(\left\vert \psi
\right\vert )$ for some smooth function $F:\left[  0,\infty\right)
\rightarrow\mathbb{R}$ and
\begin{equation}
W^{\prime}(\psi)=\frac{\partial W}{\partial\psi_{1}}+i\frac{\partial
W}{\partial\psi_{2}},\ \ \ \psi=\psi_{1}+i\psi_{2} \label{mo}%
\end{equation}
namely
\[
W^{\prime}(\psi)=F^{\prime}(\left\vert \psi\right\vert )\frac{\psi}{\left\vert
\psi\right\vert }.
\]

Equation (\ref{KG}) is the Euler-Lagrange equation relative to the Lagrangian
density
\[
\mathcal{L}=\operatorname{Re}\left(  i\partial_{t}\psi\overline{\psi}\right)
-\frac{1}{2}\left\vert \nabla\psi\right\vert ^{2}-V(x)\left\vert
\psi\right\vert ^{2}-W\left(  \psi\right)  .
\]
By the well known Noether's theorem (see e.g. \cite{gelfand}, \cite{milano})
the invariance of $\mathcal{L}$ under a one parameter Lie group gives rise to
a constant of the motion.

Since $\mathcal{L}$ is invariant under the action of the time translations the
\textit{energy}
\[
\mathcal{E(}\psi)=\frac{1}{2}\int\left[  |\nabla\psi|^{2}+V(x)\left\vert
\psi\right\vert ^{2}\right]  dx+\int W(\psi)
\]
is constant along the solutions of (\ref{KG}).

Since $W(\psi)=F(\left\vert \psi\right\vert ),$ $\mathcal{L}$ is invariant
under the $S^{1}$ action $^{{}}$
\[
\psi\rightarrow e^{i\theta}\psi,
\]
then the \textit{charge }$C(\psi),$ defined by\textit{\ }%
\[
C(\psi)=\int\left\vert \psi\right\vert ^{2},
\]
is constant along the solutions of (\ref{KG}) (see e.g. \cite{milano}).

The angular momentum, by definition, is the quantity which is preserved by
virtue of the invariance under space rotations (with respect to the origin) of
the Lagrangian. We shall consider, for simplicity, the case of three space
dimensions $N=3.$ If we assume that
\[
V(x_{1},x_{2},x_{3})=V(x_{3})
\]
namely that $V$ depends only on the third coordinate, then the Lagrangian is
invariant under the group of rotations around the axis $x_{3}.$ In this case
the third component of the momentum%

\begin{equation}
M_{3}(\psi)=\operatorname{Re}\int\left(  x_{1}\partial_{x_{2}}\psi
-x_{2}\partial_{x_{1}}\psi\right)  \;dx\text{ }%
\end{equation}
is a constant of motion. Using the polar form
\begin{equation}
\psi(t,x)=u(t,x)e^{iS(t,x)},\text{ }u\geq0 \label{polar}%
\end{equation}
$M_{3}(\psi)$ can be written as follows
\[
M_{3}(\psi)=\int\left(  x_{1}\partial_{x_{2}}S-x_{2}\partial_{x_{1}}S\right)
\ u^{2}\;dx.
\]

A solution of (\ref{KG}) is called standing wave if it has the following form:%
\begin{equation}
\psi\left(  t,x\right)  =\psi_{0}\left(  x\right)  e^{-i\omega t}%
\,,\quad\omega>0\, \label{stationary}%
\end{equation}
A \textit{vortex} is a standing wave with nonvanishing \textit{angular
momentum}.

It is immediate to check that if $\psi_{0}\left(  x\right)  $ in
(\ref{stationary}) has real values, the angular momentum $M_{3}(\psi)$ is
trivial$\mathbf{.}$ However, if $\psi_{0}\left(  x\right)  $ is allowed to
have complex values, it is possible to have $M_{3}(\psi)\neq0\mathbf{.}$ Thus,
we are led to make an ansatz of the following form:
\begin{equation}
\psi\left(  t,x\right)  =u\left(  x\right)  e^{i\left(  \ell\theta\left(
x\right)  -\omega t\right)  }\,,\quad u\left(  x\right)  \geq0,~\omega
\in\mathbb{R},\;\ell\in\mathbb{Z}-\left\{  0\right\}  \label{ansatz}%
\end{equation}
and
\[
\theta\left(  x\right)  =\operatorname{Im}\log(x_{1}+ix_{2})\in\mathbb{R}%
/2\pi\mathbb{Z};\mathbb{\;\;}x=(x_{1},x_{2},x_{3}).
\]
Moreover, we assume that $u$ has a \textit{cylindrical symmetry,}
namely\textit{\ }
\begin{equation}
u(x)=u(r,x_{3}),\ \text{where }r=\sqrt{x_{1}^{2}+x_{2}^{2}}. \label{mariella}%
\end{equation}

By this ansatz, equation (\ref{KG}) is equivalent to the system
\[
\left\{
\begin{array}
[c]{l}%
-\triangle u+\ell^{2}\left\vert \nabla\theta\right\vert ^{2}u+W^{\prime
}\left(  u\right)  +2V(x_{3})\psi=2\omega u\\
u^{2}\triangle\theta+2\nabla u\cdot\nabla\theta=0\,.
\end{array}
\right.
\]
By the definition of $\theta$ and (\ref{mariella}) we have
\[
\triangle\theta=0\,,\quad\nabla\theta\cdot\nabla u=0\,,\quad\left\vert
\nabla\theta\right\vert ^{2}=\frac{1}{r^{2}}%
\]
where the dot $\cdot$ denotes the euclidean scalar product.

So the above system reduces to find solutions, with symmetry (\ref{mariella}),
of the equation
\begin{equation}
-\triangle u+\frac{\ell^{2}}{r^{2}}u+W^{\prime}\left(  u\right)
+2V(x)u=2\omega u\qquad\text{in }\mathbb{R}^{3}. \label{elliptic eq}%
\end{equation}
Direct computations show that the energy and the third component of the
angular momentum become
\begin{align}
E_{\ell}\left(  u\right)   &  =\mathcal{E}\left(  u\left(  x\right)
e^{i\left(  \ell\theta\left(  x\right)  -\omega t\right)  }\right)
\label{energy finite}\\
&  =\int_{\mathbb{R}^{3}}\left[  \frac{1}{2}\left\vert \nabla u\right\vert
^{2}+\left(  \frac{1}{2}\frac{\ell^{2}}{r^{2}}+V(x)\right)  u^{2}+W\left(
u\right)  \right]  dx
\end{align}%
\begin{equation}
M_{3}\left(  u\left(  x\right)  e^{i\left(  \ell\theta\left(  x\right)
-\omega t\right)  }\right)  =-\ell\int_{\mathbb{R}^{3}}u^{2}dx.
\label{momentum nonzero}%
\end{equation}

We point out that $M_{3}$ in (\ref{momentum nonzero}) is nontrivial when both
$\ell$ and $u$ are not zero. Let us remark that the solutions of equation
(\ref{elliptic eq}) can be obtained as critical points of the functional
(\ref{energy finite}) on the manifold
\[
\mathfrak{M}_{c}:=\left\{  u\in X\ |\ C(u)=c\right\}
\]
where $X$ is the Hilbert space obtained by the closure of $\mathcal{D}%
(\mathbb{R}^{N})$\footnote{$\mathcal{D(\mathbb{R}}^{N}\mathcal{)}$ denotes the
space of the infinitely differentiable functions with compact support defined
in $\mathcal{\mathbb{R}}^{N}$.} with respect to the norm%
\begin{equation}
\left\Vert u\right\Vert _{X}^{2}=\int_{\mathbb{R}^{3}}\left[  \left\vert
\nabla u\right\vert ^{2}+\left(  \frac{\ell}{r^{2}}+1\right)  u^{2}\right]
dx. \label{norm}%
\end{equation}

Thus we can apply the minimization result (Theorem \ref{furbo}) stated in
section \ref{mr}. Cleary, using this approach, $2\omega$ will be the Lagrange multiplier.

\subsection{Existence of vortices\label{vort}}

Recall that $W(\psi)=F(\left\vert \psi\right\vert ).$With abuse of notation,
in the following we write $W$ instead of $F$ We make the following assumptions:

\begin{description}
\item (i) $W:\mathbb{R}_{+}\rightarrow\mathbb{R}$ is a $C^{2}$ function which
satisfies the following assumptions:
\begin{equation}
W(0)=W^{\prime}(0)=W^{\prime\prime}(0)=0 \tag{$W_0$}\label{W}%
\end{equation}%
\begin{equation}
|W^{\prime}(s)|\leq c_{1}s^{r-1}+c_{2}s^{q-1}\text{ for }q,\text{ }r\text{ in
}(2,2^{\ast}),2^{\ast}=\frac{2N}{N-2},N\geq3 \tag{$W_1$}\label{Wp}%
\end{equation}%
\begin{equation}
W(s)\geq-cs^{p},\text{ }c\geq0,\ 2<p<2+\frac{4}{N}\text{ for }s\text{ large}
\tag{$W_2$}\label{W0}%
\end{equation}%
\begin{equation}
\exists s_{0}\in\mathbb{R}^{+}\text{ such that }\frac{W(s_{0})}{s_{0}^{2}%
}<\inf\ V-\sup V \tag{$W_3$}\label{W1}%
\end{equation}

\item (ii) $V:\mathbb{R}^{N}\rightarrow\mathbb{R}$ is a continuous function
which satisfies the following assumptions:
\begin{equation}
1\leq V(x)\leq V_{0}<\infty; \tag{$V_0$}\label{V0}%
\end{equation}

\item
\begin{equation}
\forall k\in\mathbb{Z},\text{ }V(x)=V(x_{3})=V(x_{3}+k). \tag{$V_1$}\label{V1}%
\end{equation}

\end{description}

We get the following theorem:

\begin{theorem}
\label{vortices}Assume that (W$_{0}),$...,(W$_{3})$ and (\ref{V0}), ($V_{1}),
$are satisfied. Then, for any integer $\ell\neq0,$ there exist $\bar{\delta
}>0$ and a family $\psi_{\delta},$ $\delta\in(0,$ $\bar{\delta}),$ of vortices
of (\ref{KG}) with angular momentum $\left(  0,0,-\ell\int_{\mathbb{R}^{3}%
}\left\vert \psi_{\delta}\right\vert ^{2}dx\right)  .$
\end{theorem}

\begin{remark}
The conditions (\ref{W}) and (\ref{V0}) are assumed for simplicity; in fact
they can be easily weakened as follows
\[
W(0)=W^{\prime}(0)=0,\ \ W^{\prime\prime}(0)=E_{0}%
\]
and
\[
E_{1}\leq V(x)\leq V_{0}<+\infty.
\]
In fact, in the general case, it is possible to replace $W(s)$ with
\[
W_{1}(s)=W(s)-\frac{1}{2}E_{0}s^{2}%
\]
and $V(x)$ with
\[
V_{1}(x)=V(x)-E_{1}+1.
\]
In this case equation (\ref{elliptic eq}) becomes
\begin{equation}
-\triangle u+\frac{\ell^{2}}{r^{2}}u+W_{1}^{\prime}\left(  u\right)
+2V_{1}(x)u=\left(  -E_{0}-2E_{1}+2+2\omega\right)  u\qquad\text{in
}\mathbb{R}^{3}.
\end{equation}
Thus in the general case, there is only a change of the lagrange multiplier
and so the solution of the Schroedinger equation is modified only by a phase factor.
\end{remark}

By the preceding subsection we deduce that the existence of vortices of
angular momentum $\ell$ is reduced to find critical points, having the
symmetry (\ref{mariella}), of the functional $E_{\ell}$ (\ref{energy finite})
on $\mathfrak{M}_{c}$.

Now consider the action $T_{\theta}$ of the group $S^{1}$ on $u(x_{1},x_{2}, $
$x_{3})\in X,$ defined by
\begin{equation}
T_{\theta}u=u(R_{\theta}(x_{1},x_{2}),x_{3}),\ \ \theta\in\frac{\mathbb{R}%
}{2\pi\mathbb{Z}}, \label{action}%
\end{equation}
where $R_{\theta}$ denotes the rotation of an angle $\theta$ in the plane
$x_{1},$ $x_{2}.\ $We set
\[
X_{r}=\left\{  u\in X\mid u=u(r,x_{3})\right\}  ,\text{ }\mathfrak{M}_{c}%
^{r}=\mathfrak{M}_{c}\cap X_{r}.
\]
Observe that $V$ depends only on $x_{3}$, then the functional $E_{\ell}$ is
invariant under the action (\ref{action}). So by the Palais principle of
symmetric criticality \cite{Palais}, the critical points of $E_{\ell}$ on
$\mathfrak{M}_{c}^{r}$ are also critical points of $E_{\ell}$ on
$\mathfrak{M}_{c};$moreover these critical points clearly have the symmetry
(\ref{mariella}).

These observations show that the proof of Theorem \ref{vortices} is an
immediate conseguence of the following proposition:

\begin{proposition}
\label{imp}Let the assumptions of Theorem \ref{vortices} be satisfied. Then ,
for any integer $\ell,$ there exist $\bar{\delta}>0$ and a family of values of
charges $c_{\delta},$ $\delta\in(0,$ $\bar{\delta}),$ such that $E_{\ell}$
possesses a minimizer on any $\mathfrak{M}_{c_{\delta}}^{r}.$
\end{proposition}

In order to prove Proposition \ref{imp} we shall use Theorem \ref{furbo}. In
this case we have
\[
\mathbf{u}=u,\text{ }u\in X_{r}\text{ and }E(\mathbf{u)}=E_{\ell}(u),\text{
}C(\mathbf{u})=C(u)=\int u^{2}dx
\]
where $E_{\ell}(u)$ is defined in (\ref{energy finite}).

Observe first that, by (\ref{V1}), $E_{\ell}$ is invariant under the action
$T_{k}$ of the group $G=\mathbb{Z}$ on $X_{r}$ defined by%

\[
\text{ }T_{k}u(r,x_{3})=u(r,x_{3}+k),\ \ k\in\mathbb{Z}%
\]
Clearly $E_{\ell}$ and $C$ satisfy assumptions (EC-1), (EC-2).

In the following Lemmas we shall show that $E_{\ell}$ and $C$ satisfy also
(EC-3), (EC-4) and (\ref{hh}).

\begin{lemma}
\label{coercive}Let the assumptions of Theorem \ref{vortices} be satisfied.
Then $E_{\ell}\ $and $C$ satisfy the coercivity assumption (EC-3).
\end{lemma}

\textbf{Proof.} We recall a well known inequality: there exists a constant
$b_{p}>0,$ such that for any $u$ in $H^{1}(\mathbb{R}^{N})$ ($N\geq3)$
\begin{equation}
||u||_{L^{p}}\leq b_{p}||u||_{L^{2}}^{1-N\left(  \frac{1}{2}-\frac{1}%
{p}\right)  }||\nabla u||_{L^{2}}^{N\left(  \frac{1}{2}-\frac{1}{p}\right)  }.
\label{sob}%
\end{equation}
Then%
\begin{equation}
||u||_{L^{p}}^{p}\leq b_{p}||u||_{L^{2}}^{p-pN\left(  \frac{1}{2}-\frac{1}%
{p}\right)  }||\nabla u||_{L^{2}}^{pN\left(  \frac{1}{2}-\frac{1}{p}\right)
}.
\end{equation}
Since $2<p<2+\frac{4}{N},$ then $pN\left(  \frac{1}{2}-\frac{1}{p}\right)
:=q<2.$ So%
\begin{equation}
||u||_{L^{p}}^{p}\leq b_{p}||u||_{L^{2}}^{r}||\nabla u||_{L^{2}}^{q}
\label{main}%
\end{equation}
where $r=p-pN\left(  \frac{1}{2}-\frac{1}{p}\right)  =p-q>0.$

Then by H\"{o}lder inequality we have%
\begin{align*}
||u||_{L^{p}}^{p}  &  \leq b_{p}M||u||_{L^{2}}^{r}\frac{1}{M}||\nabla
u||_{L^{2}}^{q}\\
&  \leq\frac{1}{\gamma^{\prime}}\left(  b_{p}M||u||_{L^{2}}^{r}\right)
^{\gamma^{\prime}}+\frac{1}{\gamma}\left(  \frac{1}{M}||\nabla u||_{L^{2}}%
^{q}\right)  ^{\gamma}\\
&  =\frac{\left(  b_{p}M\right)  ^{\gamma^{\prime}}}{\gamma^{\prime}%
}||u||_{L^{2}}^{r\gamma^{\prime}}+\frac{1}{\gamma M^{\gamma}}||\nabla
u||_{L^{2}}^{q\gamma}.
\end{align*}
Now chose $\gamma=\frac{2}{q}$ and $M=$ $\left(  \frac{2c}{\gamma}\right)
^{1/\gamma}$ (where $c$ is the constant in assumption (W$_{2}))$ so that%
\[
||u||_{L^{p}}^{p}\leq\frac{\left(  b_{p}M\right)  ^{\gamma^{\prime}}}%
{\gamma^{\prime}}||u||_{L^{2}}^{r\gamma^{\prime}}+\frac{1}{2c}||\nabla
u||_{L^{2}}^{2}.
\]
Then
\begin{equation}
c||u||_{L^{p}}^{p}\leq a||u||_{L^{2}}^{2s}+\frac{1}{2}||\nabla u||_{L^{2}}^{2}
\label{none}%
\end{equation}
where%
\[
a=\frac{c\left(  b_{p}M\right)  ^{\gamma^{\prime}}}{\gamma^{\prime}%
};\ \ \ s=\frac{r\gamma^{\prime}}{2}.
\]
So we have, taking $N=3,$ and using $(W_{2})$ and (\ref{none})
\begin{align}
E_{\ell}(u)+aC(u)^{s}  &  =\frac{1}{2}||\nabla u||_{L^{2}}^{2}+\int
Vu^{2}+\ell^{2}\int\frac{u^{2}}{r^{2}}+\int W(u)+a||u||_{L^{2}}^{2s}%
\nonumber\\
&  \geq\frac{1}{2}||\nabla u||_{L^{2}}^{2}+\int Vu^{2}+\ell^{2}\int\frac
{u^{2}}{r^{2}}-c\int\left\vert u\right\vert ^{p}+a||u||_{L^{2}}^{2s}%
\nonumber\\
&  \geq\ell^{2}\int\frac{u^{2}}{r^{2}}+\int Vu^{2}. \label{mi}%
\end{align}
Observe that, since $p>2,$ we have $s>1.$ So (EC-3)(i) is satisfied. Moreover
it can be easily verifed that also (EC-3)(ii) is satisfied.

Now let us prove (EC-3)(iii). Let $u_{n}$ be a bounded sequence in $X_{r}$
such that $\Phi(u_{n})\rightarrow0,$ then by (\ref{mi}) we have%
\begin{equation}
\ell^{2}\int\frac{u_{n}^{2}}{r^{2}}+\int Vu_{n}^{2}\rightarrow0.
\label{partone}%
\end{equation}
So $\int u_{n}^{2}\rightarrow0$ and $\ell^{2}\int\frac{u_{n}^{2}}{r^{2}%
}\rightarrow0.$ Then, in order to show that $u_{n}$ goes to $0$ in $X_{r},$ it
remains to prove that%
\begin{equation}
||\nabla u_{n}||_{L^{2}}^{2}\rightarrow0. \label{norma}%
\end{equation}
Since $u_{n}$ is bounded in $X_{r},$ by (\ref{main}) we get
\begin{equation}
\int\left\vert u_{n}\right\vert ^{p}\rightarrow0. \label{parte}%
\end{equation}
Since $\Phi(u_{n})\rightarrow0$ and by assumption (W$_{2})$, we have%
\begin{equation}
0=\lim(E_{\ell}(u_{n})+aC(u_{n})^{s})\geq\lim\sup(\frac{1}{2}||\nabla
u_{n}||_{L^{2}}^{2}+D_{n}) \label{comp}%
\end{equation}
where%
\begin{equation}
D_{n}=\ell^{2}\int\frac{u_{n}^{2}}{r^{2}}+\int Vu_{n}^{2}+a\int u_{n}%
^{2}-c\int\left\vert u_{n}\right\vert ^{p}. \label{partissima}%
\end{equation}
By (\ref{partone}) and (\ref{parte}) we get $D_{n}\rightarrow0.$ So by
(\ref{comp}) we deduce (\ref{norma}).

$\square$

\begin{lemma}
\label{split}Let the assumptions of theorem \ref{vortices} be satisfied. Then
$E_{\ell}$ and $C$ satisfy the splitting property (EC-4).
\end{lemma}

\textbf{Proof. }Consider any sequence
\[
u_{n}=u+w_{n}\in X_{r}%
\]
where $w_{n}$ converges weakly to $0.$ We set%

\[
E_{\ell}\left(  v\right)  =A(v,v)+K(v\mathbf{)}%
\]
where%
\[
A(v,v)=\int_{\mathbb{R}^{3}}\left[  \frac{1}{2}\left\vert \nabla v\right\vert
^{2}+\left(  \frac{1}{2}\frac{\ell^{2}}{r^{2}}+V(x)\right)  v^{2}\right]
\]
and%

\[
K(v\mathbf{)=}\int W\left(  v\right)  dx.
\]
Since $C(v)=\int v^{2}$ and $A(v,v)$ are quadratic, by remark \ref{quadratic},
we have only to show that $K(v\mathbf{)}$ satisfies the splitting property.
For any measurable $A\subset\mathbb{R}^{3}$ and any $\nu\in X_{r} $, we set%
\[
K_{A}(v\mathbf{)=}\int_{A}W(v)dx.
\]
Choose $\varepsilon>0$ and $R=R(\varepsilon)>0$ such that
\begin{equation}
\left\vert K_{B_{R}^{c}}\left(  u\right)  \right\vert <\varepsilon\label{bis}%
\end{equation}
where $\ $%
\[
B_{R}^{c}=\mathbb{R}^{N}-B_{R}\text{ and }B_{R}=\left\{  x\in\mathbb{R}%
^{N}:\left\vert x\right\vert <R\right\}  .
\]
Since $w_{n}\rightharpoonup0$ weakly in $H^{1}\left(  \mathbb{R}^{3}\right)
$, by usual compactness arguments, we have that
\begin{equation}
K_{B_{R}}\left(  w_{n}\right)  \rightarrow0\text{ and }K_{B_{R}}\left(
u+w_{n}\right)  \rightarrow K_{B_{R}}\left(  u\right)  . \label{bibis}%
\end{equation}
Then, by (\ref{bis}) and (\ref{bibis}), we have%

\begin{align}
&  \underset{n\rightarrow\infty}{\lim}\left\vert K\left(  u+w_{n}\right)
-K\left(  u\right)  -K\left(  w_{n}\right)  \right\vert \nonumber\\
&  =\ \underset{n\rightarrow\infty}{\lim}\left\vert K_{B_{R}^{c}}\left(
u+w_{n}\right)  +K_{B_{R}}\left(  u+w_{n}\right)  -K_{B_{R}^{c}}\left(
u\right)  -K_{B_{R}}\left(  u\right)  -K_{B_{R}^{c}}\left(  w_{n}\right)
-K_{B_{R}}\left(  w_{n}\right)  \right\vert \nonumber\\
&  \mathbb{=}\ \underset{n\rightarrow\infty}{\lim}\left\vert K_{B_{R}^{c}%
}\left(  u+w_{n}\right)  -K_{B_{R}^{c}}\left(  u\right)  -K_{B_{R}^{c}}\left(
w_{n}\right)  \right\vert \nonumber\\
&  \leq\ \underset{n\rightarrow\infty}{\lim}\left\vert K_{B_{R}^{c}}\left(
u+w_{n}\right)  -K_{B_{R}^{c}}\left(  w_{n}\right)  \right\vert +\varepsilon.
\label{pepe}%
\end{align}

Now, by the intermediate value theorem, there exists $\zeta_{n}\in(0,1)$ such
that for $z_{n}=$ $\zeta_{n}u+\left(  1-\zeta_{n}\right)  w_{n}$, we have that%

\[
\left\vert K_{B_{R}^{c}}\left(  u+w_{n}\right)  -K_{B_{R}^{c}}\left(
w_{n}\right)  \right\vert =\left\vert \left\langle K_{B_{R}^{c}}^{\prime
}\left(  z_{n}\right)  ,u\right\rangle \right\vert
\]
\begin{align}
&  \leq\int_{B_{R}^{c}}\left\vert W^{\prime}(z_{n})u\right\vert \leq(\text{by
(\ref{Wp})})\\
&  \leq\int_{B_{R}^{c}}c_{1}\left\vert z_{n}\right\vert ^{r-1}\left\vert
u\right\vert +c_{2}\left\vert z_{n}\right\vert ^{q-1}\left\vert u\right\vert
\\
&  \leq c_{1}\left\Vert z_{n}\right\Vert _{L^{r}(B_{R}^{c})}^{r-1}\left\Vert
u\right\Vert _{L^{r}(B_{R}^{c})}+c_{2}\left\Vert z_{n}\right\Vert
_{L^{q}(B_{R}^{c})}^{q-1}\left\Vert u\right\Vert _{_{L^{q}(B_{R}^{c})}}\\
&  (\text{if }R\text{ is large enough)}\\
&  \leq c_{3}\left(  \left\Vert z_{n}\right\Vert _{L^{r}(B_{R}^{c})}%
^{r-1}+\left\Vert z_{n}\right\Vert _{L^{q}(B_{R}^{c})}^{q-1}\right)
\varepsilon
\end{align}
So we have
\begin{equation}
\left\vert K_{B_{R}^{c}}\left(  u+w_{n}\right)  -K_{B_{R}^{c}}\left(
w_{n}\right)  \right\vert \leq c_{3}\left(  \left\Vert z_{n}\right\Vert
_{L^{r}(B_{R}^{c})}^{r-1}+\left\Vert z_{n}^{{}}\right\Vert _{L^{q}(B_{R}^{c}%
)}^{q-1}\right)  \varepsilon\label{media}%
\end{equation}

Since $z_{n}$ is bounded in $H^{1}\left(  \mathbb{R}^{3}\right)  ,$ the
sequences $\left\Vert z_{n}\right\Vert _{L^{r}(B_{R}^{c})}^{r-1}$ and
$\left\Vert z_{n}^{{}}\right\Vert _{L^{q}(B_{R}^{c})}^{q-1}$are bounded. Then,
by (\ref{pepe}) and (\ref{media}), we easily get%
\begin{equation}
\underset{n\rightarrow\infty}{\lim}\left\vert K\left(  u+w_{n}\right)
-K\left(  u\right)  -K\left(  w_{n}\right)  \right\vert \leq\varepsilon
+M\cdot\varepsilon. \label{inno}%
\end{equation}

where $M$ is a suitable constant.

Since $\varepsilon$ is arbitrary, from (\ref{inno}) we get
\[
\underset{n\rightarrow\infty}{\lim}\left\vert K\left(  u+w_{n}\right)
-K\left(  u\right)  -K\left(  w_{n}\right)  \right\vert =0
\]

$\square$

Now in order to prove that assumption (\ref{hh}) is satisfied some work is
necessary. Set
\[
\Lambda(u)=\frac{E_{\ell}(u)}{C(u)}.
\]
First of all we have:

\begin{lemma}
\label{preparatorio}If the assumptions of Theorem \ref{vortices} are
satisfied, then for $6>t>2,$ we have
\[
\underset{u\in X_{r},\left\Vert u\right\Vert _{L^{t}}\rightarrow0}{\lim\inf
}\Lambda(u)=\underset{u\in X_{r}\left\Vert u\right\Vert _{Lt}=1}{\inf}%
\frac{\frac{1}{2}\int\left(  \left\vert \nabla u\right\vert ^{2}+\frac
{\ell^{2}u^{2}}{r^{2}}\right)  dx+\int Vu^{2}}{\int u^{2}}.
\]

\end{lemma}

\textbf{Proof. }Clearly%

\[
\underset{u\in X_{r},\left\Vert u\right\Vert _{L^{t}}\rightarrow0}{\lim\inf
}\Lambda(u)=\ \underset{u\in X_{r},\left\Vert u\right\Vert _{L^{t}%
}=1,\varepsilon\rightarrow0}{\lim\inf}\frac{E(\varepsilon u)}{C(\varepsilon
u)}%
\]%
\[
=\underset{u\in X_{r},\left\Vert u\right\Vert _{L^{t}}=1,}{\inf}\left(
\frac{\frac{1}{2}\int\left(  \left\vert \nabla u\right\vert ^{2}+\frac
{\ell^{2}u^{2}}{r^{2}}\right)  dx+\int Vu^{2}}{\int u^{2}}\right)
+\underset{u\in X_{r},\left\Vert u\right\Vert _{L^{t}}=1,\varepsilon
\rightarrow0}{\lim\inf}\frac{\int W(\varepsilon u)}{\varepsilon^{2}\int u^{2}%
}.
\]

So the proof of Lemma will be achieved if we show that%
\begin{equation}
\underset{u\in X_{r},\left\Vert u\right\Vert _{L^{t}}=1,\varepsilon
\rightarrow0}{\lim\inf}\frac{\int W(\varepsilon u)}{\varepsilon^{2}\int u^{2}%
}=0. \label{resto}%
\end{equation}
By assumptions $\left(  W_{1}\right)  $ and $\left(  W_{2}\right)  $ we have%
\begin{equation}
-cs^{p}\leq W(s)\leq\bar{c}(s^{q}+s^{r}) \label{zero}%
\end{equation}
where $c,\bar{c}$ are positive constants and $q,r$ in $(2,6).$

Then by (\ref{zero}) we have%
\begin{equation}
-cA\varepsilon^{p-2}\leq\underset{\left\Vert u\right\Vert _{L^{t}}=1}{\inf
}\frac{\int W(\varepsilon u)}{\varepsilon^{2}\int u^{2}}\leq\bar
{c}B(\varepsilon^{q-2}+\varepsilon^{r-2}) \label{uno}%
\end{equation}
where%
\[
A=\underset{u\in X_{r}\text{ }\left\Vert u\right\Vert _{L^{t}}=1}{\inf}%
\frac{\int\left\vert u\right\vert ^{p}}{\int u^{2}},\text{ }B=\underset{u\in
X_{r}\text{ }\left\Vert u\right\Vert _{L^{t}}=1}{\inf}\frac{\int\left(
\left\vert u\right\vert ^{q}+\left\vert u\right\vert ^{r}\right)  }{\int
u^{2}}.
\]
By (\ref{uno}) we easily get (\ref{resto}).

$\square$

Now consider the following action $T_{k}$ of the group $G=\mathbb{Z}$ on
$X_{r}:$
\[
\text{for all }u\in X_{r},\text{ }k\in\mathbb{Z}\text{ }T_{k}u(x_{1}%
,x_{2},x_{3})=u(x_{1},x_{2},x_{3}+k)
\]

The following proposition holds

\begin{lemma}
\label{nonvanishing}\textit{\textbf{\ }}If\textit{\textbf{\ }}$2<t<6,$ the
norm $\left\Vert u\right\Vert _{L^{t}}$ satisfies the property (\ref{seminorm}%
), namely%
\[
\left\{  u_{n}\ \text{is a vanishing sequence}\right\}  \Rightarrow\left\Vert
u_{n}\right\Vert _{L^{t}}\rightarrow0.
\]

\end{lemma}

\textbf{Proof. }Let $u_{n}\ $be a G-vanishing sequence in $X_{r}$ and, arguing
by contradiction, assume that $\left\Vert u_{n}\right\Vert _{L^{t}}$ does not
converge to $0.$ Then, up to a subsequence,%
\begin{equation}
\left\Vert u_{n}\right\Vert _{L^{t}}\geq a>0. \label{dis}%
\end{equation}
Since $u_{n}$ is bounded in $X_{r},$ we have that for a suitable constant
$M>0$
\begin{equation}
\left\Vert u_{n}\right\Vert _{H^{1}}^{2}\leq M. \label{bounded}%
\end{equation}

Now we set
\[
Q_{i}=\left\{  \left(  x_{1},x_{2},x_{3}\right)  :i\leq x_{3}<i+1\right\}
\text{, }i\text{ integer.}%
\]
Clearly
\[
\mathbb{R}^{3}=\underset{}{%
{\displaystyle\bigcup\limits_{i\in\mathbb{Z}}}
}Q_{i}%
\]
Let $C$ denote the constant for the Sobolev embedding $H^{1}\left(
Q_{i}\right)  \subset L^{t}\left(  Q_{i}\right)  ,$then, by (\ref{dis}) and
(\ref{bounded}), we get the following
\begin{align*}
0  &  <a^{t}\leq\int\left\vert u_{n}\right\vert ^{t}=\sum_{i}\int_{Q_{i}%
}\left\vert u_{n}\right\vert ^{t}=\sum_{i}\left\Vert u_{n}\right\Vert
_{L^{t}\left(  Q_{i}\right)  }^{t-2}\left\Vert u_{n}\right\Vert _{L^{t}\left(
Q_{i}\right)  }^{2}\\
&  \leq\ \left(  \underset{i}{\sup}\left\Vert u_{n}\right\Vert _{L^{t}\left(
Q_{i}\right)  }^{t-2}\right)  \cdot\sum_{i}\left\Vert u_{n}\right\Vert
_{L^{t}\left(  Q_{i}\right)  }^{2}\\
&  \leq\ C\left(  \underset{i}{\sup}\left\Vert u_{n}\right\Vert _{L^{t}\left(
Q_{i}\right)  }^{t-2}\right)  \cdot\sum_{i}\left\Vert u_{n}\right\Vert
_{H^{1}\left(  Q_{i}\right)  }^{2}\\
&  =C\left(  \underset{i}{\sup}\left\Vert u_{n}\right\Vert _{L^{t}\left(
Q_{i}\right)  }^{t-2}\right)  \left\Vert \mathbf{u}_{n}\right\Vert _{H^{1}%
}^{2}\leq CM\left(  \underset{i}{\sup}\left\Vert u_{n}\right\Vert
_{L^{t}\left(  Q_{i}\right)  }^{t-2}\right)  .
\end{align*}
Then%
\[
\left(  \underset{i}{\sup}\left\Vert u_{n}\right\Vert _{L^{t}\left(
Q_{i}\right)  }\right)  \geq\left(  \frac{a^{t}}{CM}\right)  ^{1/(t-2)}.
\]

So, for any $n,$ there exists an integer $i_{n}$ such that
\begin{equation}
\left\Vert u_{n}\right\Vert _{L^{t}\left(  Q_{i_{n}}\right)  }\geq\alpha>0.
\end{equation}
Then
\[
\left\Vert T_{i_{n}}u_{n}\right\Vert _{L^{t}(Q_{0})}=\left\Vert u_{n}%
\right\Vert _{L^{t}(Q_{i_{n}})}\geq\alpha>0.
\]
Since $u_{n}$ and then $T_{_{i_{n}}}u_{n}$ is bounded in $X_{r}$, we have,
passing eventually to a subsequence, that
\[
T_{_{i_{n}}}u_{n}\rightharpoonup u_{0}\text{ weakly in }X_{r}.
\]
Clearly, if we show that $u_{0}\neq0,$ we get a contradiction with the
assumption that $u_{n}$ is nonvanishing.

Now, let $\varphi=\varphi\left(  x_{3}\right)  $ be a nonnegative, $C^{\infty
}$-function whose value is $1$ for $0<x_{3}<1$ and $0$ for $\left\vert
x_{3}\right\vert >2.$ Then the sequence $\varphi T_{i_{n}}u_{n}$ is bounded in
$H_{0}^{1}(\mathbb{R}^{2}\times(-2,2)),$ moreover $\varphi T_{i_{n}}u_{n}$ is
invariant under the action (\ref{action}). Then, using the compactness result
proved in \cite{el}, we have
\[
\varphi T_{i_{n}}u_{n}\rightarrow\chi\text{ strongly in }L^{t}(\mathbb{R}%
^{2}\times(-2,2)).
\]
On the other hand
\begin{equation}
\varphi T_{i_{n}}u_{n}\rightarrow\varphi u_{0}\text{ }a.e\text{.}
\label{punto}%
\end{equation}

Then
\begin{equation}
\varphi T_{i_{n}}u_{n}\rightarrow\varphi u_{0}\text{ strongly in }%
L^{t}(\mathbb{R}^{2}\times(-2,2)). \label{virgola}%
\end{equation}

Moreover
\begin{equation}
\left\Vert \varphi T_{i_{n}}u_{n}\right\Vert _{L^{3}\left(  \mathbb{R}%
^{2}\times(-2,2)\right)  }\geq\left\Vert \varphi T_{i_{n}}u_{n}\right\Vert
_{L^{t}\left(  Q_{0}\right)  }=\left\Vert u_{n}\right\Vert _{L^{t}\left(
Q_{i_{n}}\right)  }\geq\alpha>0. \label{pp}%
\end{equation}

Then$\ $by (\ref{virgola}) and (\ref{pp})
\[
\left\Vert \varphi u_{0}\right\Vert _{L^{3}\left(  \mathbb{R}^{2}%
\times(-2,2)\right)  }\geq\alpha>0.
\]
Thus we have that $u_{0}\neq0.$

$\square$

Finally it remains to show that also assumption (\ref{hh}) is satisfied.

\begin{lemma}
\label{seguente}Let the assumptions of Theorem \ref{vortices} be satisfied,
then
\[
\underset{u\in X_{r}}{\inf}\Lambda(u)<\Lambda_{0}.
\]

\end{lemma}

\textbf{Proof. }We shall show that there exists $u\in X_{r}$ such that
$\Lambda(u)<\Lambda_{0}.$ The construction of such $u$ needs some work since
we require that $u$ belongs to $X_{r},$ namely we require that $u$ is
invariant under the $S^{1}$action (\ref{action}) and it is $0$ near the
$x_{3}$ axis, so that $\int\frac{u^{2}}{r^{2}}$ converges.

For $0<\mu<\lambda$ we set:
\[
T_{\lambda,\mu}=\left\{  \left(  r,x_{3}\right)  :(r-\lambda)^{2}+x_{3}{}%
^{2}\leq\mu^{2}\right\}
\]
and, for $\lambda>2,$ we consider a smooth function $u_{\lambda}$ with
cylindrical symmetry such that
\begin{equation}
u_{\lambda}(r,x_{3})=\left\{
\begin{array}
[c]{cc}%
s_{0} & if\;\;\left(  r,x_{3}\right)  \in T_{\lambda,\lambda/2}\\
& \\
0 & if\;\;\left(  r,x_{3}\right)  \notin T_{\lambda,\lambda/2+1}%
\end{array}
\right.  \label{ni}%
\end{equation}
where $s_{0}$ is such that $\frac{W(s_{0})}{s_{0}^{2}}<\inf V-\sup V$ (see
(W$_{3}$)). Moreover we may assume that
\begin{equation}
\left\vert \nabla u_{\lambda}\left(  r,x_{3}\right)  \right\vert
\leq2\;\text{for }\left(  r,x_{3}\right)  \in T_{\lambda,\lambda
/2+1}\backslash T_{\lambda,\lambda/2}. \label{nu}%
\end{equation}

We have
\begin{equation}
\Lambda(u_{\lambda})=\frac{\int\left[  \left\vert \nabla u_{\lambda
}\right\vert ^{2}+\frac{\ell^{2}u_{\lambda}^{2}}{r^{2}}+2Vu^{2}\right]
dx}{2\int u_{\lambda}^{2}}+\frac{\int W(u_{\lambda})dx}{\int u_{\lambda}^{2}}.
\label{lan}%
\end{equation}

By (\ref{nu}) and (\ref{ni}) a direct computation shows that
\begin{equation}
\int\left\vert \nabla u_{\lambda}\right\vert ^{2}\leq4meas(T_{\lambda
,\lambda/2+1}\backslash T_{\lambda,\lambda/2})\leq c_{1}\lambda^{2} \label{al}%
\end{equation}%
\begin{equation}
\int\frac{u_{\lambda}^{2}}{r^{2}}\leq\frac{c_{2}}{\lambda^{2}}meas(T_{\lambda
,\lambda/2+1})\leq c_{3}\lambda\label{all}%
\end{equation}%
\begin{equation}
\int u_{\lambda}^{2}\geq c_{4}meas(T_{\lambda,\lambda/2})\geq c_{5}\lambda^{3}
\label{alll}%
\end{equation}

where $c_{1},..,$ $c_{5}$ are positive constants. So that
\begin{equation}
\frac{\int\left[  \left\vert \nabla u_{\lambda}\right\vert ^{2}+\frac{\ell
^{2}u_{\lambda}^{2}}{r^{2}}+2Vu_{\lambda}^{2}\right]  dx}{2\int u_{\lambda
}^{2}}\leq\sup V+O\left(  \frac{1}{\lambda}\right)  . \label{pin}%
\end{equation}
Now
\[
\int_{T_{\lambda,\lambda/2+1}\backslash T_{\lambda,\lambda/2}}\left\vert
W(u_{\lambda})\right\vert \leq c_{6}meas(T_{\lambda,\lambda/2+1}\backslash
T_{\lambda,\lambda/2})\leq c_{6}\lambda^{2}%
\]
Then
\begin{align*}
\int W(u_{\lambda})dx  &  =W(s_{0})meas(T_{\lambda,\lambda/2})+\int
_{T_{\lambda,\lambda/2+1}\backslash T_{\lambda,\lambda/2}}W(u_{\lambda})\\
&  \leq W(s_{0})meas(T_{\lambda,\lambda/2})+\int_{T_{\lambda,\lambda
/2+1}\backslash T_{\lambda,\lambda/2}}\left\vert W(u_{\lambda})\right\vert \\
&  \leq W(s_{0})meas(T_{\lambda,\lambda/2})+c_{6}\lambda^{2}.
\end{align*}
So%
\begin{align}
\frac{\int W(u_{\lambda})dx}{\int u_{\lambda}^{2}}  &  \leq\frac
{W(s_{0})meas(T_{\lambda,\lambda/2})+c_{6}\lambda^{2}}{\int_{{}}u_{\lambda
}^{2}}\leq\nonumber\\
&  \frac{W(s_{0})meas(T_{\lambda,\lambda/2})}{\int_{{}}u_{\lambda}^{2}}%
+\frac{c_{6}\lambda^{2}}{\int_{{}}u_{\lambda}^{2}}. \label{ala}%
\end{align}

Now, since $W(s_{0})<0,$ we have
\begin{equation}
\frac{W(s_{0})meas(T_{\lambda,\lambda/2})}{\int_{{}}u_{\lambda}^{2}}\leq
\frac{W(s_{0})meas(T_{\lambda,\lambda/2})}{s_{0}^{2}meas(T_{\lambda
,\lambda/2+1})}=\frac{W(s_{0})}{s_{0}^{2}}\left(  \frac{\lambda}{\lambda
+2}\right)  ^{2}. \label{ali}%
\end{equation}
Then by (\ref{alll}), (\ref{ala}) and (\ref{ali}) we have%
\begin{equation}
\frac{\int W(u_{\lambda})dx}{\int u_{\lambda}^{2}}\leq\frac{W(s_{0})}%
{s_{0}^{2}}\left(  \frac{\lambda}{\lambda+2}\right)  ^{2}+\frac{c_{7}}%
{\lambda}. \label{alto}%
\end{equation}
By (\ref{lan}), (\ref{pin}) and (\ref{alto}) we get%
\begin{align*}
\Lambda(u_{\lambda})  &  \leq\sup V+\frac{W(s_{0})}{s_{0}^{2}}\left(
\frac{\lambda}{\lambda+2}\right)  ^{2}+O\left(  \frac{1}{\lambda}\right)  .\\
&  \sup V-\inf V+\inf V+\frac{W(s_{0})}{s_{0}^{2}}\left(  \frac{\lambda
}{\lambda+2}\right)  ^{2}+O\left(  \frac{1}{\lambda}\right)  .
\end{align*}
By lemma \ref{preparatorio} we have
\[
\underset{\mathbf{u}\in X_{r},\left\Vert \mathbf{u}\right\Vert _{L^{t}%
}\rightarrow0}{\lim\inf}\Lambda(\mathbf{u})\geq\inf V
\]
then
\begin{equation}
\Lambda(u_{\lambda})\leq\sup V-\inf V+\underset{\mathbf{u}\in X_{r},\left\Vert
\mathbf{u}\right\Vert _{L^{t}}\rightarrow0}{\lim\inf}\Lambda(\mathbf{u}%
)+\frac{W(s_{0})}{s_{0}^{2}}\left(  \frac{\lambda}{\lambda+2}\right)
^{2}+O\left(  \frac{1}{\lambda}\right)  . \label{fere}%
\end{equation}

By assumption (W$_{3})$ for $\lambda$ large we have
\begin{equation}
\sup V-\inf V+\frac{W(s_{0})}{s_{0}^{2}}\left(  \frac{\lambda}{\lambda
+2}\right)  ^{2}+O\left(  \frac{1}{\lambda}\right)  <0. \label{fer}%
\end{equation}
By (\ref{fere}) and (\ref{fer}) we get that for $\lambda$ large%
\begin{equation}
\Lambda(u_{\lambda})<\underset{\mathbf{u}\in X_{r},\left\Vert u\right\Vert
_{L^{t}}\rightarrow0}{\lim\inf}\Lambda(u). \label{quasi}%
\end{equation}

On the other hand, since by Lemma \ref{nonvanishing} $\left\Vert {}\right\Vert
_{L^{t}}$ satisfies the property (\ref{seminorm}), we have by (\ref{aaaa})
that%
\begin{equation}
\underset{u\in X_{r},\left\Vert u\right\Vert _{L^{t}}\rightarrow0}{\lim\inf
}\Lambda(u)\leq\Lambda_{0}. \label{quasiquasi}%
\end{equation}

Clearly (\ref{quasi}) and (\ref{quasiquasi}) imply that assumption (\ref{hh})
is satisfied.

$\square$

\section{Existence of hylomorphic solitons\label{ehs}}

\subsection{An abstract definition of solitary waves and solitons\label{be}}

Solitons are particular \textit{states} of a dynamical system described by one
or more partial differential equations. Thus we assume that the states of this
system are described by one or more \textit{fields} which mathematically are
represented by functions
\begin{equation}
\mathbf{u}:\mathbb{R}^{N}\rightarrow V \label{lilla}%
\end{equation}
where $V$ is a vector space with norm $\left\vert \ \cdot\ \right\vert _{V}$
and which is called the internal parameters space. We assume the system to be
deterministic; this means that it can be described as a dynamical system
$\left(  X,\gamma\right)  $ where $X$ is the set of the states and
$\gamma:\mathbb{R}\times X\rightarrow X$ is the time evolution map. If
$\mathbf{u}_{0}(x)\in X,$ the evolution of the system will be described by the
function
\begin{equation}
\mathbf{u}\left(  t,x\right)  :=\gamma_{t}\mathbf{u}_{0}(x). \label{flusso}%
\end{equation}
We assume that the states of $X$ have "finite energy" so that they decay at
$\infty$ sufficiently fast. Roughly speaking, the solitons are "bump"
solutions characterized by some form of stability.

To define them at this level of abstractness, we need to recall some well
known notions in the theory of dynamical systems.

\begin{definition}
A set $\Gamma\subset X$ is called \textit{invariant} if $\forall\mathbf{u}%
\in\Gamma,\forall t\in\mathbb{R},\ \gamma_{t}\mathbf{u}\in\Gamma.$
\end{definition}

\begin{definition}
Let $\left(  X,d\right)  $ be a metric space and let $\left(  X,\gamma\right)
$ be a dynamical system. An invariant set $\Gamma\subset X$ is called stable,
if $\forall\varepsilon>0,$ $\exists\delta>0,\;\forall\mathbf{u}\in X$,
\[
d(\mathbf{u},\Gamma)\leq\delta,
\]
implies that
\[
\forall t\geq0,\text{ }d(\gamma_{t}\mathbf{u,}\Gamma)\leq\varepsilon.
\]

\end{definition}

Let $G$ be a subgroup of $(\mathbb{R}^{N},+)$ and consider the following
action $T_{z}$ of $G$ on $X:$ for all $z\in G$ and $\mathbf{u}\in X$
\[
T_{z}\mathbf{u}(x)=\mathbf{u}(x+z).
\]
$\ $

Now we are ready to give the definition of soliton:

\begin{definition}
\label{ds} A state $\mathbf{u}(x)\in X$ is called soliton if there is an
invariant set $\Gamma$ such that

\begin{itemize}
\item (i) $\forall t,\ \gamma_{t}\mathbf{u}(x)\in\Gamma,$

\item (ii) $\Gamma$ is stable,

\item (iii) $\Gamma$ is $G$-compact (Def. \ref{gcompatto})
\end{itemize}
\end{definition}

\begin{remark}
The above definition needs some explanation. For simplicity, we assume that
$\Gamma$ is a manifold (actually, in many concrete models, this is the generic
case). Then (iii) implies that $\Gamma$ is finite dimensional. Since $\Gamma$
is invariant, $\mathbf{u}_{0}\in\Gamma\Rightarrow\gamma_{t}\mathbf{u}_{0}%
\in\Gamma$ for every time. Thus, since $\Gamma$ is finite dimensional, the
evolution of $\mathbf{u}_{0}$ is described by a finite number of parameters$.$
The dynamical system $\left(  \Gamma,\gamma\right)  $\ behaves as a point in a
finite dimensional phase space. By the stability of $\Gamma$, a small
perturbation of $\mathbf{u}_{0}$ remains close to $\Gamma.$ However, in this
case, its evolution depends on an infinite number of parameters. Thus, this
system appears as a finite dimensional system with a small perturbation.
\end{remark}

\begin{remark}
The type of stability described above is called \textbf{orbital stability} in
the literature relative to the nonlinear Schr\"{o}dinger and Klein-Gordon
equations\textbf{.}
\end{remark}

\subsection{An existence result for hylomorphic solitons}

We now assume that the dynamical system $\left(  X,\gamma\right)  $ has two
constants of motion. These constants can be considered as functionals on $X.$
One of them will be called energy and it will be denoted by $E$; the other
will be called hylenic charge and it will be denoted by $C$.

At this level of abstractness, the names energy and hylenic charge are
conventional but $E$ and $C$ satisfy different assumptions; see assumption
(EC-3) in section \ref{af}. In our applications to PDE's,\ $E$ will be the
usual energy. The name hylenic charge has been introduced in \cite{milano},
\cite{BBBM} and \cite{hylo}.

The presence of $E$ and $C$ allows to give the following definition of
hylomorphic soliton.

\begin{definition}
\label{tdc}A soliton $\mathbf{u}_{0}\in X$ is called hylomorphic if $\Gamma$
(as in Def. \ref{ds}) has the following structure%
\[
\Gamma=\Gamma\left(  e_{0},c_{0}\right)  =\left\{  \mathbf{u}\in
X\ |\ E(\mathbf{u})=e_{0},\ C(\mathbf{u})=c_{0}\right\}
\]
where%
\begin{equation}
e_{0}=\min\left\{  E(\mathbf{u})\ |\ C(\mathbf{u})=c_{0}\right\}  \label{min}%
\end{equation}
for some $c_{0}\in\mathbb{R}$.
\end{definition}

Notice that, by (\ref{min}), we have that a hylomorphic soliton $\mathbf{u}%
_{0}$ satisfies the following nonlinear eigenvalue problem:%
\[
E^{\prime}(\mathbf{u}_{0})=\lambda C^{\prime}(\mathbf{u}_{0}).
\]

Clearly, for a given $c_{0}$ the minimum $e_{0}$ in (\ref{min}) might not
exist; moreover, even if the minimum exists, it is possible that $\Gamma
\ $does not satisfy (ii) or (iii) of def. \ref{ds}.

The following theorem holds

\begin{theorem}
\label{astra1} Assume that the dynamical system $(X,\gamma)$ satisfies
(EC-1),...,(EC-4) and (\ref{hh}). Moreover assume that%
\[
E\text{ and }C\text{ are two constants of motion.}%
\]
\newline Then there exists $\bar{\delta}>0$ such that the dynamical system
$(X,\gamma)$ admits a family $\mathbf{u}_{\delta}$ ($\delta\in\left(
0,\bar{\delta}\right)  $) of hylomorphic solitons.
\end{theorem}

The proof of this theorem will be given in the next section.

\subsection{A stability result}

In order to prove Theorem \ref{astra1} it is sufficient to show that the
minimizers in Th. \ref{furbo} provide solitons, so we have to prove that the
set $\Gamma_{c_{\delta}}$ is stable. To do this, we need the (well known)
Liapunov theorem in following form:

\begin{theorem}
\label{propV}Let $\Gamma$ be an invariant set and assume that there exists a
differentiable function $V$ (called a Liapunov function) such that

\begin{itemize}
\item (a) $V(\mathbf{u})\geq0$ and\ $V(\mathbf{u})=0\Leftrightarrow u\in
\Gamma$

\item (b) $\partial_{t}V(\gamma_{t}\left(  \mathbf{u}\right)  )\leq0$

\item (c) $V(\mathbf{u}_{n})\rightarrow0\Leftrightarrow d(\mathbf{u}%
_{n},\Gamma)\rightarrow0.$
\end{itemize}

\noindent Then $\Gamma$ is stable.
\end{theorem}

\textbf{Proof. }For completeness, we give a proof of this well known result.
Arguing by contradiction, assume that $\Gamma,$ satisfying the assumptions of
Th. \ref{propV}, is not stable. Then there exists $\varepsilon>0$ and
sequences $\mathbf{u}_{n}\in X$ and $t_{n}>0$ such that
\begin{equation}
d(\mathbf{u}_{n},\Gamma)\rightarrow0\text{ and }d(\gamma_{t_{n}}\left(
\mathbf{u}_{n}\right)  ,\Gamma)>\varepsilon. \label{bingo}%
\end{equation}
Then we have%
\[
d(\mathbf{u}_{n},\Gamma)\rightarrow0\Longrightarrow V(\mathbf{u}%
_{n})\rightarrow0\Longrightarrow V(\gamma_{t_{n}}\left(  \mathbf{u}%
_{n}\right)  )\rightarrow0\Longrightarrow d(\gamma_{t_{n}}\left(
\mathbf{u}_{n}\right)  ,\Gamma)\rightarrow0
\]
where the first and the third implications are consequence of property (c).
The second implication follows from property (b). Clearly, this fact
contradicts (\ref{bingo}).

$\square$

\begin{theorem}
\label{astraco}Assume (EC-1) and (EC-2). For $\mathbf{u}\in X$ and
$e_{0},c_{0}\in$ $\mathbb{R}$, we set%
\begin{equation}
V\left(  \mathbf{u}\right)  =\left(  E\left(  \mathbf{u}\right)
-e_{0}\right)  ^{2}+\left(  C\left(  \mathbf{u}\right)  -c_{0}\right)  ^{2}.
\end{equation}
If $V$ is$\ G$-compact (see Def. \ref{gcompattoa}) and
\begin{equation}
\Gamma=\left\{  \mathbf{u}\in X:E(\mathbf{u})=e_{0},\ C(\mathbf{u}%
)=c_{0}\right\}  \neq\varnothing, \label{gamma}%
\end{equation}
then every $\mathbf{u}\in\Gamma$ is a soliton.
\end{theorem}

\textbf{Proof}: We have to prove that $\Gamma$ in (\ref{gamma}) satisfies
(i),(ii) and (iii) of Def. \ref{ds}. The property (iii), namely the fact that
$\Gamma$ is G-compact, is a trivial consequence of the fact that $\Gamma$ is
the set of minimizers of a G-compact functional $V$ (see definitions
\ref{gcompatto} and \ref{gcompattoa}). The invariance property (i) is clearly
satisfied since $E$ and $C$ are constants of the motion. It remains to prove
(ii), namely that $\Gamma$ is stable. To this end we shall use Th.
\ref{propV}. So we need to show that $V(\mathbf{u})$ satisfies (a), (b) and
(c). Statements (a) and (b) are trivial. Now we prove (c). First we show the
implication $\Rightarrow.$ Let $\mathbf{u}_{n}$ be a sequence such that
$V(\mathbf{u}_{n})\rightarrow0.$ By contradiction we assume that
$d(\mathbf{u}_{n},\Gamma)\nrightarrow0,$ namely that there is a subsequence
$\mathbf{u}_{n}^{^{\prime}}$ such that
\begin{equation}
d(\mathbf{u}_{n}^{\prime},\Gamma)\geq a>0. \label{kaka}%
\end{equation}
Since $V(\mathbf{u}_{n})\rightarrow0$ also $V(\mathbf{u}_{n}^{\prime
})\rightarrow0,$ and, since $V$ is $G$ compact, there exists a sequence
$g_{n}$ in $G$ such that, for a subsequence $\mathbf{u}_{n}^{\prime\prime}$,
we have $g_{n}\mathbf{u}_{n}^{\prime\prime}\rightarrow\mathbf{u}_{0}.$ Then
\[
d(\mathbf{u}_{n}^{\prime\prime},\Gamma)=d(g_{n}\mathbf{u}_{n}^{\prime\prime
},\Gamma)\leq d(g_{n}\mathbf{u}_{n}^{\prime\prime},\mathbf{u}_{0}%
)\rightarrow0
\]
and this contradicts (\ref{kaka}).

Now we prove the other implication $\Leftarrow.$ Let $\mathbf{u}_{n}$ be a
sequence such that $d(\mathbf{u}_{n},\Gamma)\rightarrow0,$ then there exists
$\mathbf{v}_{n}\in\Gamma$ s.t.
\begin{equation}
d(\mathbf{u}_{n},\Gamma)\geq d(\mathbf{u}_{n},\mathbf{v}_{n})-\frac{1}{n}.
\label{triplo}%
\end{equation}

Since $V$ is G-compact, also $\Gamma$ is G-compact; so, for a suitable
sequence $g_{n}$, we have $g_{n}\mathbf{v}_{n}\rightarrow\mathbf{\bar{w}}%
\in\Gamma.$ We get the conclusion if we show that $V(\mathbf{u}_{n}%
)\rightarrow0.$ We have by (\ref{triplo}), that $d(\mathbf{u}_{n}%
,\mathbf{v}_{n})\rightarrow0$ and hence $d(g_{n}\mathbf{u}_{n},g_{n}%
\mathbf{v}_{n})\rightarrow0$ and so, since $g_{n}\mathbf{v}_{n}\rightarrow
\mathbf{\bar{w},}$ we have $g_{n}\mathbf{u}_{n}\rightarrow\mathbf{\bar{w}}%
\in\Gamma.$ Therefore, by the continuity of $V$ and since $\mathbf{\bar{w}}%
\in\Gamma,$ we have $V\left(  g_{n}\mathbf{u}_{n}\right)  \rightarrow V\left(
\mathbf{\bar{w}}\right)  =0$ and we can conclude that $V\left(  \mathbf{u}%
_{n}\right)  \rightarrow0.$

$\square$

\bigskip

In the cases in which we are interested, $X$ is an infinite dimensional
manifold; then if you choose generic $e_{0}$ and $c_{0},$ $V$ is not
$G$-compact since the set $\Gamma=\left\{  \mathbf{u}\in X:E(\mathbf{u}%
)=e_{0},\ C(\mathbf{u})=c_{0}\right\}  $ has codimension 2. However, Th.
(\ref{furbo}) allows to determine $e_{0}$ and $c_{0}$ in such a way that $V$
is $G $-compact and hence to prove the existence of solitons by using Theorem
\ref{astraco}.

\bigskip

\textbf{Proof of Th. \ref{astra1}.} In order to prove Th. \ref{astra1}, we
will use Th. \ref{astraco} with $e_{0}=e_{\delta}$ and $c_{0}=c_{\delta}$
where $e_{\delta}$ and $c_{\delta}$ are given by Th. \ref{furbo}.

We set
\begin{equation}
V\left(  \mathbf{u}\right)  =\left(  E\left(  \mathbf{u}\right)  -e_{\delta
}\right)  ^{2}+\left(  C\left(  \mathbf{u}\right)  -c_{\delta}\right)  ^{2}.
\label{liap}%
\end{equation}
We show that $V$ is $G$-compact: let $\mathbf{w}_{n}$ be a minimizing sequence
for $V,$ then $V\left(  \mathbf{w}_{n}\right)  \rightarrow0$ and consequently
$E\left(  \mathbf{w}_{n}\right)  \rightarrow e_{\delta}$ and $C\left(
\mathbf{w}_{n}\right)  \rightarrow c_{\delta}$. Let $J_{\delta}$ be as in
Theorem \ref{furbo}. Now, since
\[
\inf J_{\delta}=\frac{e_{\delta}}{c_{\delta}}+\delta\left[  e_{\delta
}+ac_{\delta}^{s}\right]  ,
\]
we have that $\mathbf{w}_{n}$ is a minimizing sequence also for $J_{\delta}. $
Then, since $J_{\delta}$ is $G$-compact, we get
\begin{equation}
\mathbf{w}_{n}\ \text{is}\ G\text{-compact}. \label{paracula}%
\end{equation}
So we conclude that $V$ is $G$-compact and hence the conclusion follows by
using Theorem \ref{astraco}.

$\square$

\section{ Existence of solitons for NSE with periodic
potential\label{periodic}}

In this section we shall study the existence of hylomorphic solitons on
lattice for the Schr\"{o}dinger equation (\ref{KG}) in $\mathbb{R}^{N}.$

The existence of solitons for (\ref{KG}) is an old problem and there are many
results in the case $V=0$ (\cite{CL82}, \cite{bus03}, \cite{BBGM} and the
references in \cite{milano}).

Here we assume that $V$ is a lattice potential, namely we assume that the
potential $V$ satisfies the periodicity condition.%

\begin{equation}
V(x)=V(x+Az)\text{ for all }x\in\mathbb{R}^{N}\text{ and }z\in\mathbb{Z}^{N}
\tag{$V1^\prime$}\label{V1'}%
\end{equation}
where $A$ is a $N\times N$ invertible matrix.

Here we look for solitons and do not require they to be vortices, so the
energy corresponds to the expression (\ref{energy finite}) with $\ell=0,$
namely%
\begin{equation}
E(u)=E_{0}\left(  u\right)  =\int_{\mathbb{R}^{3}}\left[  \frac{1}%
{2}\left\vert \nabla u\right\vert ^{2}+V(x)u^{2}+W\left(  u\right)  \right]
dx,\text{ }u\in X. \label{m}%
\end{equation}

As before the charge is
\[
C(u)=\int u^{2}dx.
\]
In this case $X$ is the ordinary $H^{1}(\mathbb{R}^{N})$ Sobolev space. We
shall consider the following action of the group $G=\mathbb{Z}^{N}$ on $X:$%
\begin{equation}
\text{for all }z\in\mathbb{Z}^{N}\text{ and }u\in X:\text{ }T_{z}u(x)=u(x+Az).
\label{group}%
\end{equation}
Clearly the charge $C$ is $G$-invariant and, since $V$ satisfies (\ref{V1'}),
also the energy $E$ is invariant under this group action. The following
Theorem holds:

\begin{theorem}
\label{solitoni}Let $W$and $V$ satisfy assumptions (W$_{0}),...,(W_{3})$ and
(V$_{0}),$(V$_{1}^{\prime}).$ Then there exists $\bar{\delta}>0$ such that the
dynamical system described by the Schr\"{o}dinger equation (\ref{KG}) has a
family $u_{\delta}$ ($\delta\in\left(  0,\bar{\delta}\right)  )$ of
hylomorphic solitons.
\end{theorem}

The proof of this theorem is based on the abstract theorem \ref{astra1}. In
this case the energy is given by (\ref{m}). We need to show that assumptions
$(W_{0}),...,(W_{3})$ and $(V_{0}),\ (V_{1}^{\prime})$ permit to show that
assumptions (EC-1), ...,(EC-4) and (\ref{hh}) of theorem \ref{astra1} are satisfied.

(EC-1) and (EC-2) are trivially verified. The proof of the other assumptions
follows the same lines of the proof of Th. \ref{vortices} as we can see in the
following lemmas:

\begin{lemma}
\label{coercive2} $E\ $and $C$ satisfy the coercivity assumption (EC-3).
\end{lemma}

\textbf{Proof.} The proof is the same of that of lemma \ref{coercive} with
$\ell=0$.

$\square$

\begin{lemma}
\label{split2} $E$ and $C$ satisfy the splitting property (EC-4).
\end{lemma}

\textbf{Proof.} The proof is the same of that of lemma \ref{split} with
$\ell=0$.

$\square$

\begin{lemma}
\label{nonvanishing2}\textit{\textbf{\ }}If\textit{\textbf{\ }}$2<t<\frac
{2N}{N-2},$ $N\geq3,$ the norm $\left\Vert u\right\Vert _{L^{t}}$ satisfies
the property (\ref{seminorm}), namely%
\[
\left\{  \mathbf{u}_{n}\ \text{is a vanishing sequence}\right\}
\Rightarrow\left\Vert u_{n}\right\Vert _{L^{t}}\rightarrow0.
\]

\end{lemma}

\textbf{Proof. } We set for $j\in\mathbb{Z}^{N}$
\[
Q_{j}=A\left(  j+Q^{0}\right)  =\left\{  Aj+Aq:q\in Q^{0}\right\}
\]
where $Q^{0}$ is now the cube defined as follows
\[
Q^{0}=\left\{  \left(  x_{1},..,x_{n}\right)  \in\mathbb{R}^{N}:0\leq
x_{i}<1\right\}  \text{.}%
\]
Now let $x\in\mathbb{R}^{N}$ and set $y=A^{-1}(x).$ Clearly there exist $q\in
Q^{0}$ and $j\in\mathbb{Z}^{N}$ such that $y=j+q.$ So
\[
x=Ay=A(j+q)\in Q_{j}.
\]
Then we conclude that
\[
\mathbb{R}^{N}=%
{\displaystyle\bigcup\limits_{j}}
Q_{j}.
\]
Let $u_{n}\ $be a bounded sequence in $H^{1}\left(  \mathbb{R}^{N}\right)  $
such that, up to a subsequence, $\left\Vert u_{n}\right\Vert _{L^{t}}\geq
a>0.$ We need to show that $u_{n}$ is non vanishing. Then, if $C$ is the
constant for the Sobolev embedding $H^{1}\left(  Q_{j}\right)  \subset
L^{t}\left(  Q_{j}\right)  $ and $\left\Vert u_{n}\right\Vert _{H^{1}}^{2}\leq
M,$ we have
\begin{align*}
0  &  <a^{t}\leq\int\left\vert u_{n}\right\vert ^{t}=\sum_{j}\int_{Q_{j}%
}\left\vert u_{n}\right\vert ^{t}=\sum_{j}\left\Vert u_{n}\right\Vert
_{L^{t}\left(  Q_{j}\right)  }^{t-2}\left\Vert u_{n}\right\Vert _{L^{t}\left(
Q_{j}\right)  }^{2}\\
&  \leq\ \left(  \underset{j}{\sup}\left\Vert u_{n}\right\Vert _{L^{t}\left(
Q_{j}\right)  }^{t-2}\right)  \cdot\sum_{j}\left\Vert u_{n}\right\Vert
_{L^{t}\left(  Q_{j}\right)  }^{2}\\
&  \leq\ C\left(  \underset{j}{\sup}\left\Vert u_{n}\right\Vert _{L^{t}\left(
Q_{j}\right)  }^{t-2}\right)  \cdot\sum_{j}\left\Vert u_{n}\right\Vert
_{H^{1}\left(  Q_{j}\right)  }^{2}\\
&  =C\left(  \underset{j}{\sup}\left\Vert u_{n}\right\Vert _{L^{t}\left(
Q_{j}\right)  }^{t-2}\right)  \left\Vert u_{n}\right\Vert _{H^{1}}^{2}\leq
CM\left(  \underset{j}{\sup}\left\Vert u_{n}\right\Vert _{L^{t}\left(
Q_{j}\right)  }^{t-2}\right)  .
\end{align*}
Then%
\[
\left(  \underset{j}{\sup}\left\Vert u_{n}\right\Vert _{L^{t}\left(
Q_{j}\right)  }\right)  \geq\left(  \frac{a^{t}}{CM}\right)  ^{1/(t-2)}%
\]

Then, for any $n,$ there exists $j_{n}\in\mathbb{Z}^{N}$ such that
\begin{equation}
\left\Vert u_{n}\right\Vert _{L^{t}\left(  Q_{j_{n}}\right)  }\geq\alpha>0.
\label{caca}%
\end{equation}
Then, if we set $Q=AQ^{0},$we easily have
\begin{equation}
\left\Vert T_{j_{n}}u_{n}\right\Vert _{L^{t}(Q)}=\left\Vert u_{n}\right\Vert
_{L^{t}(Q_{j_{n}})}\geq\alpha>0. \label{chicco}%
\end{equation}

Since $u_{n}$ is bounded, also $T_{j_{n}}u_{n}$ is bounded in $H^{1}%
(\mathbb{R}^{N}).$ Then we have, up to a subsequence, that $T_{j_{n}}%
u_{n}\rightharpoonup u_{0}$ weakly in $H^{1}(\mathbb{R}^{N})$ and hence
strongly in $L^{t}(Q)$. By (\ref{chicco}), $u_{0}\neq0.$

$\square$

\begin{lemma}
\label{seguente2}Assumption (\ref{hh}) is satisfied namely%
\[
\underset{u\in H^{1}(R^{N})}{\inf}\Lambda(u)<\Lambda_{0}%
\]

\end{lemma}

\textbf{Proof. } This lemma is analogous to lemma \ref{seguente}, however in
this case the proof is easier: since $X=H^{1}(\mathbb{R}^{N}),$ we need only
to construct a function $u\in H^{1}(\mathbb{R}^{N})$ such that $\Lambda
(u)<\Lambda_{0}.$

Such a function can be constructed as follows. Set%
\[
u_{R}=\left\{
\begin{array}
[c]{cc}%
s_{0} & if\;\;|x|<R\\
0 & if\;\;|x|>R+1\\
\frac{|x|}{R}s_{0}-(\left\vert x\right\vert -R)\frac{R+1}{R}s_{0} &
if\;\;R<|x|<R+1
\end{array}
.\right.
\]
Then
\[
\int\left\vert \nabla u_{R}\right\vert ^{2}dx=O(R^{N-1}),\int\left\vert
u_{R}\right\vert ^{2}dx=O(R^{N}),
\]

so that%

\begin{equation}
\frac{\int\left[  \left\vert \nabla u_{R}\right\vert ^{2}+2Vu_{R}^{2}\right]
dx}{2\int u_{R}^{2}}\leq\sup V+O\left(  \frac{1}{R}\right)  . \label{pinco}%
\end{equation}
Moreover
\[
\int W(u_{R})dx=W(s_{0})meas(B_{R})+\int_{B_{R+1}\backslash B_{R}}W(u_{R}).
\]
So%
\[
\frac{\int W(u_{R})dx}{\int u_{R}^{2}}\leq\frac{W(s_{0})meas(B_{R}%
)+c_{1}R^{N-1}}{\int_{{}}u_{\lambda}^{2}}\leq(\text{ since }W(s_{0})<0)
\]%
\begin{equation}
\leq\frac{W(s_{0})meas(B_{R})}{s_{0}^{2}meas(B_{R+1})}+\frac{c_{2}R^{N-1}%
}{R^{N}}=\frac{W(s_{0})}{s_{0}^{2}}\left(  \frac{R}{R+1}\right)  ^{N}%
+\frac{c_{2}}{R}. \label{palle}%
\end{equation}

Then, by (\ref{pinco}) e (\ref{palle}) we get%
\[
\Lambda(u_{R})\leq\sup V+\frac{W(s_{0})}{s_{0}^{2}}\left(  \frac{R}%
{R+1}\right)  ^{N}+O\left(  \frac{1}{R}\right)  .
\]
By lemma \ref{preparatorio} we have $\inf V\leq\underset{\left\Vert
u\right\Vert _{L^{t}}\rightarrow0}{\text{ }\lim\inf}\Lambda(u),$ then
\begin{equation}
\Lambda(u_{R})\leq\underset{\left\Vert u\right\Vert _{L^{t}}\rightarrow
0}{\text{ }\lim\inf}\Lambda(u)+\sup V-\inf V+\frac{W(s_{0})}{s_{0}^{2}}\left(
\frac{R}{R+1}\right)  ^{N}+O\left(  \frac{1}{R}\right)  . \label{ma}%
\end{equation}

On the other hand, since by Lemma \ref{nonvanishing} the $L^{t}$ norm
satisfies the property (\ref{seminorm}), we have by (\ref{aaaa}) that%
\begin{equation}
\underset{\left\Vert u\right\Vert _{L^{t}}\rightarrow0}{\lim\inf}%
\Lambda(u)\leq\Lambda_{0}. \label{me}%
\end{equation}

Clearly (\ref{ma}), (\ref{me}) and assumption (W$_{3})$ imply that for $R$
large we have%

\[
\Lambda(u_{R})<\Lambda_{0}.
\]
Then assumption (\ref{hh}) is satisfied.

$\square$

\textbf{Proof of Th. \ref{solitoni}}. The proof follows from Th. \ref{astra1}
and Lemmas \ref{coercive2}, \ref{split2} and \ref{seguente2}.

$\square$

\section{Existence of Solitons for the nonlinear Klein-Gordon
equation\label{kleingordon}}

In this section we shall apply the abstract theorem \ref{astra1} to the
existence of hylomorphic solitons in $\mathbb{R}^{N}$ for the nonlinear
Klein-Gordon equation (NKG). There are well known results on the existence of
stable solutions for (NKG) (\cite{shatah}, \cite{gss87}) and more recently the
existence of hylomorphic solitons for (NKG) has been studied in \cite{BBBM}
and \ref{hylo}.

More exactly, we consider the equation%
\begin{equation}
\square\psi+W^{\prime}(\psi)=0 \tag{NKG}\label{NKG}%
\end{equation}
where $\square=\partial_{t}^{2}-\Delta$,$\;\psi:\mathbb{R\times R}%
^{N}\rightarrow\mathbb{C}$ ($N\geq3$) , $W:\mathbb{C}\rightarrow\mathbb{R} $
and $W^{\prime}$ are as in (\ref{mo}) (see the beginning of section
\ref{beginning}). Assume that%
\begin{equation}
W(s)=\frac{1}{2}\ m^{2}s^{2}+N(s),\ \ s\geq0,\ m\neq0 \label{h0}%
\end{equation}
where $\ N(s)=o(s^{2}).$

We make the following assumptions on $W$:

\begin{itemize}
\item (NKG-i) \textbf{(Positivity}) $W(s)\geq0$ for $s\geq0$

\item (NKG-ii) \textbf{(Hylomorphy}) $\exists s_{0}\in\mathbb{R}%
^{+}\mathbb{\ }$such that $W(s_{0})<\frac{1}{2}m^{2}s_{0}^{2}$

\item (NKG-iii)\textbf{(Growth condition}) there are constants $c_{1}%
,c_{2}>0,$ $2<r,q<2N/(N-2)$ such that for any $s>0:$%
\[
|N^{\prime}(s)|\ \leq c_{1}s^{r-1}+c_{2}s^{q-1}.
\]

\end{itemize}

We shall assume that the initial value problem is well posed for (NKG). Eq.
(\ref{NKG}) is the Euler-Lagrange equation of the action functional
\begin{equation}
\mathcal{S}(\psi)=\int\left(  \frac{1}{2}\left\vert \partial_{t}%
\psi\right\vert ^{2}-\frac{1}{2}|\nabla\psi|^{2}-W(\psi)\right)  dxdt.
\label{az}%
\end{equation}
The energy and the charge take the following form:
\begin{equation}
E(\psi)=\int\left[  \frac{1}{2}\left\vert \partial_{t}\psi\right\vert
^{2}+\frac{1}{2}\left\vert \nabla\psi\right\vert ^{2}+W(\psi)\right]  dx
\label{energy}%
\end{equation}%
\begin{equation}
C(\psi)=-\operatorname{Re}\int i\partial_{t}\psi\overline{\psi}\;dx.
\label{im1}%
\end{equation}
(the sign "minus"in front of the integral is a useful convention).

\subsection{The NKG as a dynamical system}

We set%
\[
X=H^{1}(\mathbb{R}^{N},\mathbb{C})\times L^{2}(\mathbb{R}^{N},\mathbb{C})
\]
and we will denote the generic element of $X$ by $\mathbf{u}=(\psi\left(
x\right)  ,\hat{\psi}\left(  x\right)  );$ then, by the well posedness
assumption, for every $\mathbf{u}\in X,$ there is a unique solution
$\psi(t,x)$ of (\ref{NKG}) such that%
\begin{align}
\psi(0,x)  &  =\psi\left(  x\right) \label{in}\\
\partial_{t}\psi(0,x)  &  =\hat{\psi}\left(  x\right)  .\nonumber
\end{align}

Using this notation, we can write equation (\ref{NKG}) in Hamiltonian form:
\begin{align}
\partial_{t}\psi &  =\hat{\psi}\label{ham1}\\
\partial_{t}\hat{\psi}  &  =\Delta\psi-W^{\prime}(\psi). \label{ham2}%
\end{align}
The time evolution map $\gamma:\mathbb{R}\times X\rightarrow X$ is defined by
\[
\gamma_{t}\mathbf{u}_{0}(x)=\mathbf{u}\left(  t,x\right)
\]
where $\mathbf{u}_{0}(x)=(\psi\left(  x\right)  ,\hat{\psi}\left(  x\right)
)\in X$ and $\mathbf{u}\left(  t,x\right)  =(\psi\left(  t,x\right)
,\hat{\psi}\left(  t,x\right)  )$ is the unique solution of (\ref{ham1}) and
(\ref{ham2}) satisfying the initial conditions (\ref{in}). The energy and the
charge, as functionals defined in $X,$ become%

\begin{equation}
E(\mathbf{u})=\int\left[  \frac{1}{2}\left\vert \hat{\psi}\right\vert
^{2}+\frac{1}{2}\left\vert \nabla\psi\right\vert ^{2}+W(\psi)\right]  dx
\end{equation}%
\begin{equation}
C(\mathbf{u})=-\operatorname{Re}\int i\hat{\psi}\overline{\psi}\;dx.
\end{equation}

\subsection{Existence results for NKG}

The following Theorem holds:

\begin{theorem}
\label{PTS}Assume that $W$ satisfies (NKG-i),...,(NKG-iii). Then there exists
$\bar{\delta}>0$ such that the dynamical system described by the equation
(NKG) has a family $\mathbf{u}_{\delta}$ ($\delta\in\left(  0,\bar{\delta
}\right)  )$ of hylomorphic solitons.
\end{theorem}

The proof of this theorem is based on the abstract theorem \ref{astra1}. In
this case the energy $E$ and the hylenic charge $C$ have the form
(\ref{energy}) and (\ref{im1}) respectively.

Assumption (EC-1) is clearly satisfied. $E$ and $C$ are invariant under
translations, so assumption (EC-2) is satisfied with respect to the action
$T_{z}$ of the group $G=\mathbb{R}^{N}$ where%
\[
T_{z}\mathbf{u}(x)=\mathbf{u}(x+z),\ \ \ z\in\mathbb{R}^{N}.
\]
It can be seen that the coercitivity assumption (EC-3) is satisfied with
$a=0.$ Arguing as in lemma \ref{split} (replacing $W$ by $N$) it can be shown
that also (EC-4) is satisfied, namely that $E$ and $C$ satisfy the splitting property.

It remains to prove (\ref{hh}). First of all we set:%
\[
\left\Vert \mathbf{u}\right\Vert _{\sharp}=\left\Vert (\psi,\hat{\psi
})\right\Vert _{\sharp}=\max\left(  \left\Vert \psi\right\Vert _{L^{r}%
},\left\Vert \psi\right\Vert _{L^{q}}\right)
\]
where $r,$ $q$ are introduced in (NKG-iii). With some abuse of notation we
shall write $\max\left(  \left\Vert \psi\right\Vert _{L^{r}},\left\Vert
\psi\right\Vert _{L^{q}}\right)  =\left\Vert \psi\right\Vert _{\sharp}$.\ 

\bigskip

\begin{lemma}
\label{nonvanishing3}\textit{\textbf{\ }}The norm $\left\Vert \mathbf{u}%
\right\Vert _{\sharp}$ satisfies the property (\ref{seminorm}), namely%
\[
\left\{  \mathbf{u}_{n}\ \text{is a vanishing sequence}\right\}
\Rightarrow\left\Vert \psi_{n}\right\Vert _{\sharp}\rightarrow0.
\]

\end{lemma}

\textbf{Proof. } Let $\psi_{n}\ $be a bounded sequence in $H^{1}\left(
\mathbb{R}^{N}\right)  $ such that, up to a subsequence, $\left\Vert \psi
_{n}\right\Vert _{\sharp}\geq a>0.$ We need to show that $\psi_{n}$ is non
vanishing. May be taking a subsequence, we have that at least one of the
following holds:

\begin{itemize}
\item (i) $\left\Vert \psi_{n}\right\Vert _{\sharp}=\left\Vert \psi
_{n}\right\Vert _{L^{r}}$

\item (ii) $\left\Vert \psi_{n}\right\Vert _{\sharp}=\left\Vert \psi
_{n}\right\Vert _{L^{q}}$
\end{itemize}

Suppose that (i) holds. Then, we argue as il lemma \ref{nonvanishing2}. If
(ii) holds, we argue in the same way replacing $r$ with $q.$

$\square$

Now we set%

\[
\Lambda_{\sharp}=\ \underset{\left\Vert \mathbf{u}\right\Vert _{\sharp
}\rightarrow0}{\lim\inf}\Lambda(\mathbf{u})=\ \underset{\varepsilon
\rightarrow0}{\lim}\ \inf\left\{  \Lambda(\psi,\hat{\psi})\ |\ \hat{\psi}\in
L^{2};\psi\in H^{1};\ \left\Vert \psi\right\Vert _{\sharp}<\varepsilon
\right\}  .
\]
By remark \ref{valutazione}, we have that $\Lambda_{0}\geq\Lambda_{\sharp};$
so let us evaluate $\Lambda_{\sharp}$.

\begin{lemma}
\label{preparatorio2}If $W$ satisfies assumption (NKG-iii), then the following
inequality holds
\[
\Lambda_{\sharp}\geq m.
\]

\end{lemma}

\textbf{Proof.} By (NKG-iii) we have%
\begin{align*}
\left\vert \int N(\left\vert \psi\right\vert )dx\right\vert  &  \leq k_{1}%
\int\left\vert \psi\right\vert ^{r}+k_{2}\int\left\vert \psi\right\vert
^{q}\ \\
&  \leq k_{1}\left\Vert \psi\right\Vert _{\sharp}^{r}+k_{2}\left\Vert
\psi\right\Vert _{\sharp}^{q}.
\end{align*}
If we assume that $\left\Vert \psi\right\Vert _{\sharp}=1$,%
\[
\left\vert \int N(\left\vert \varepsilon\psi\right\vert )dx\right\vert \leq
k_{1}\varepsilon^{r}+k_{2}\varepsilon^{q}.
\]
Now, choose $s\in\left(  2,\min(r,q)\right)  .$ Thus, if $\varepsilon>0$ is
small enough, since $r,q>s>2$, we have%
\begin{align*}
&  \varepsilon^{s}\int\left(  \left\vert \nabla\psi\right\vert ^{2}%
+m^{2}\left\vert \psi\right\vert ^{2}\right)  dx-\left\vert \int N(\left\vert
\varepsilon\psi\right\vert )dx\right\vert \\
&  \geq\varepsilon^{s}k_{3}\left\Vert \psi\right\Vert _{\sharp}^{2}%
-k_{1}\varepsilon^{r}-k_{2}\varepsilon^{q}=k_{3}\varepsilon^{s}-k_{1}%
\varepsilon^{r}-k_{2}\varepsilon^{q}\geq0
\end{align*}
and so%
\begin{equation}
\left\vert \int N(\varepsilon\left\vert \psi\right\vert )dx\right\vert
\leq\varepsilon^{s}\int\left(  \left\vert \nabla\psi\right\vert ^{2}%
+m^{2}\left\vert \psi\right\vert ^{2}\right)  dx. \label{cinzia2}%
\end{equation}

Now we clearly have%

\begin{equation}
\Lambda_{\sharp}=\ \underset{\varepsilon\rightarrow0}{\lim}\ \inf\left\{
\Lambda(\varepsilon\psi,\hat{\psi})\ |\ \hat{\psi}\in L^{2},\psi\in
H^{1},\left\Vert \psi\right\Vert _{\sharp}=1\right\}  . \label{lin}%
\end{equation}

Let us estimate $\Lambda(\varepsilon\psi,\hat{\psi})$ using (\ref{cinzia2}):%

\begin{align*}
\Lambda(\varepsilon\psi,\hat{\psi})  &  =\frac{\frac{1}{2}\int\left(
\left\vert \hat{\psi}\right\vert ^{2}+\left\vert \nabla\varepsilon
\psi\right\vert ^{2}+m^{2}\left\vert \varepsilon\psi\right\vert ^{2}\right)
dx+\int N(\left\vert \varepsilon\psi\right\vert )dx}{\left\vert
\operatorname{Re}\int i\hat{\psi}\varepsilon\overline{\psi}\;dx\right\vert }\\
&  \geq\frac{\frac{1}{2}\int\left\vert \hat{\psi}\right\vert ^{2}+\left(
\frac{\varepsilon^{2}}{2}-\varepsilon^{s}\right)  \int\left(  \left\vert
\nabla\psi\right\vert ^{2}+m^{2}\left\vert \psi\right\vert ^{2}\right)
}{\varepsilon\left\vert \operatorname{Re}\int i\hat{\psi}\overline{\psi
}\;dx\right\vert }\\
&  \geq\frac{\frac{1}{2}\int\left\vert \hat{\psi}\right\vert ^{2}%
+\frac{\varepsilon^{2}}{2}\left(  1-2\varepsilon^{s-2}\right)  m^{2}%
\int\left\vert \psi\right\vert ^{2}}{\varepsilon\left(  \int\left\vert
\hat{\psi}\right\vert ^{2}dx\right)  ^{1/2}\left(  \int\left\vert
\psi\right\vert ^{2}dx\right)  ^{1/2}}\\
&  \geq\frac{\left(  \int\left\vert \hat{\psi}\right\vert ^{2}dx\right)
^{1/2}\cdot\varepsilon m\sqrt{1-2\varepsilon^{s-2}}\left(  \int\left\vert
\psi\right\vert ^{2}dx\right)  ^{1/2}}{\varepsilon\left(  \int\left\vert
\hat{\psi}\right\vert ^{2}dx\right)  ^{1/2}\left(  \int\left\vert
\psi\right\vert ^{2}dx\right)  ^{1/2}}=m\sqrt{1-2\varepsilon^{s-2}}.
\end{align*}
Then
\begin{equation}
\lim\Lambda(\varepsilon\psi,\hat{\psi})\geq m. \label{nuovo}%
\end{equation}

So the conclusion follows by (\ref{lin}) and (\ref{nuovo}).

$\square$

Next we will show that the hylomorphy assumption (\ref{hh}) is satisfied.

\begin{lemma}
\label{ve}Assume that $W$ satisfies (NKG-i),...,(NKG-iii), then%
\[
\underset{\mathbf{u}\in X}{\inf}\Lambda(\mathbf{u})<\Lambda_{0}.
\]

\end{lemma}

\textbf{Proof}. Let $R>0;$ set
\begin{equation}
u_{R}=\left\{
\begin{array}
[c]{cc}%
s_{0} & if\;\;|x|<R\\
0 & if\;\;|x|>R+1\\
\frac{|x|}{R}s_{0}-(\left\vert x\right\vert -R)\frac{R+1}{R}s_{0} &
if\;\;R<|x|<R+1
\end{array}
\right.  . \label{inff}%
\end{equation}
By (NKG-ii) there exists $0<\beta<m$ such that
\begin{equation}
W(s_{0})\leq\frac{\beta^{2}s_{0}^{2}}{2}. \label{medio}%
\end{equation}
We set $\psi=u_{R},\ $and\ $\hat{\psi}=\beta u_{R}.$

Then
\begin{align*}
\underset{\mathbf{u}\in X}{\inf}\Lambda(\mathbf{u})  &  =\underset{\psi
,\hat{\psi}}{\ \inf}\ \frac{\int\left(  \frac{1}{2}\left\vert \hat{\psi
}\right\vert ^{2}+\frac{1}{2}\left\vert \nabla\psi\right\vert ^{2}%
+W(\psi)\right)  dx}{\left\vert \operatorname{Re}\int i\hat{\psi}%
\overline{\psi}\;dx\right\vert }\\
&  \leq\frac{\int\left(  \frac{1}{2}\beta^{2}\left\vert u_{R}\right\vert
^{2}+\frac{1}{2}\left\vert \nabla u_{R}\right\vert ^{2}+W(u_{R})\right)
dx}{\beta\int\left\vert u_{R}\right\vert ^{2}\;dx}\\
&  \leq\ \frac{\int_{\left\vert x\right\vert <R}\left(  \frac{1}{2}\beta
^{2}\left\vert u_{R}\right\vert ^{2}+W(u_{R})\right)  dx}{\beta\int
_{\left\vert x\right\vert <R}\left\vert u_{R}\right\vert ^{2}\;dx}\\
&  +\frac{\int_{R<\left\vert x\right\vert <R+1}\left(  \frac{1}{2}\beta
^{2}\left\vert u_{R}\right\vert ^{2}+\frac{1}{2}\left\vert \nabla
u_{R}\right\vert ^{2}+W(u_{R})\right)  dx}{\beta\int_{\left\vert x\right\vert
<R}\left\vert u_{R}\right\vert ^{2}dx}\\
&  =\frac{1}{2}\beta+\frac{\int_{\left\vert x\right\vert <R}W(s_{0})dx}%
{\beta\int_{\left\vert x\right\vert <R}\left\vert s_{0}\right\vert ^{2}%
dx}+O\left(  \frac{1}{R}\right)  .
\end{align*}

Then, by (\ref{medio}), we have%
\[
\underset{\mathbf{u}\in X}{\inf}\Lambda(\mathbf{u})\leq\frac{1}{2}\beta
+\frac{\int_{\left\vert x\right\vert <R}\frac{1}{2}s_{0}^{2}\beta^{2}}%
{\beta\int_{\left\vert x\right\vert <R}\left\vert s_{0}\right\vert ^{2}%
dx}+O\left(  \frac{1}{R}\right)  =\beta+O\left(  \frac{1}{R}\right)  .
\]
So
\begin{equation}
\underset{\mathbf{u}\in X}{\inf}\Lambda(\mathbf{u})\leq\beta. \label{pri}%
\end{equation}

Then, since $\beta<m,$ by remark \ref{valutazione}, lemma \ref{preparatorio2}
and (\ref{pri}), we have that
\[
\Lambda_{0}\geq\Lambda_{\sharp}\geq m>\beta\geq\ \underset{\mathbf{u}\in
X}{\inf}\Lambda(\mathbf{u})
\]
and so the conclusion easily follows.

$\square$

\textbf{Proof of Th. \ref{PTS}. }The assumptions (EC-1),...,(EC-4) and
(\ref{hh}) are satisfied, then the proof follows by using Th. \ref{astra1}.

$\square$

We conclude this section with the following theorem which gives some more
information on the structure of the solitons:

\begin{theorem}
\label{pipi}Let $\mathbf{u}$ be a hylomorphic soliton relative to the equation
(NKG) with initial data $\mathbf{u}_{0}(x)=(\psi_{0}(x),\hat{\psi}_{0}(x))$
$\in X$. Then there exists $\omega\in\mathbb{R}$ such that $\psi_{0}$
satisfies the equation
\begin{equation}
-\Delta\psi_{0}+W^{\prime}(\psi_{0})=\omega^{2}\psi_{0}, \label{ella}%
\end{equation}%
\[
\hat{\psi}_{0}=-i\omega\psi_{0}%
\]
and
\begin{equation}
\gamma_{t}\mathbf{u}_{0}(x)=\left[
\begin{array}
[c]{c}%
\psi_{0}(x)e^{-i\omega t}\\
-i\omega\psi_{0}(x)e^{-i\omega t}%
\end{array}
\right]  . \label{evol}%
\end{equation}

\end{theorem}

\textbf{Proof.} Since $\mathbf{u}$ is a hylomorphic soliton it is a critical
point of $E$ constrained on the manifold $\mathfrak{M}_{c}=\left\{
\mathbf{u}\in X:C(\mathbf{u})=c\right\}  .$

Clearly
\begin{equation}
E^{\prime}(u_{0})=-\omega C^{\prime}(u_{0}) \label{become}%
\end{equation}

where $-\omega$ is a Lagrange multiplier. We now compute the derivatives
$E^{\prime}(u_{0}),C^{\prime}(u_{0}).$

For all $(v_{0},v_{1})\in X=H^{1}(\mathbb{R}^{N},\mathbb{C})\times
L^{2}(\mathbb{R}^{N},\mathbb{C}),$ we have%
\[
E^{\prime}(u_{0})\left[
\begin{array}
[c]{c}%
v_{0}\\
v_{1}%
\end{array}
\right]  =\operatorname{Re}\int\left[  \hat{\psi}_{0}\overline{v_{1}}%
+\nabla\psi_{0}\overline{\nabla v_{0}}+W^{\prime}(x,\psi_{0})\overline{v_{0}%
}\right]  dx
\]%
\begin{align*}
C^{\prime}(u_{0})\left[
\begin{array}
[c]{c}%
v_{0}\\
v_{1}%
\end{array}
\right]   &  =-\operatorname{Re}\int\left(  i\hat{\psi}_{0}\overline{v_{0}%
}+iv_{1}\overline{\psi_{0}}\right)  \;dx\\
&  =-\operatorname{Re}\int\left(  i\hat{\psi}_{0}\overline{v_{0}}%
+\overline{iv_{1}\overline{\psi_{0}}}\right)  \;dx\\
&  =-\operatorname{Re}\int\left(  i\hat{\psi}_{0}\overline{v_{0}}-i\psi
_{0}\overline{v_{1}}\right)  \;dx.
\end{align*}
Then (\ref{become}) can be written as follows:
\begin{align*}
\operatorname{Re}\int\left[  \nabla\psi_{0}\overline{\nabla v_{0}}+W^{\prime
}(x,\psi_{0})\overline{v_{0}}\right]  dx  &  =\omega\operatorname{Re}\int
i\hat{\psi}_{0}\overline{v_{0}}\;dx\\
\operatorname{Re}\int\hat{\psi}_{0}\overline{v_{1}}\ dx  &  =-\omega
\operatorname{Re}\int i\psi_{0}\overline{v_{1}}\;dx.
\end{align*}
Then%
\begin{align}
-\Delta\psi_{0}+W^{\prime}(x,\psi_{0})  &  =i\omega\hat{\psi}_{0}\nonumber\\
\hat{\psi}_{0}  &  =-i\omega\psi_{0}. \label{lui}%
\end{align}

So we get (\ref{ella}). From (\ref{ella}) and (\ref{lui}), we easily verify
that (\ref{evol}) solves (\ref{ham1}), (\ref{ham2}).

$\square$

\end{document}